\documentclass{article}

\usepackage{amssymb}

\title{Model-completion of scaled lattices\thanks{This 
paper is an revised version of \cite{prep-latt-2004}, 
submitted to APAL the 5th of April 2005 and eventually rejected 
the 24th of March 2006.\hfill\break
Keywords: model-theory, quantifier elimination, scaled lattice, Heyting algebra, p-adic.\hfill\break
MSC classes: 03C10, 06D20, 06D99.}}

\author{Luck Darni\`ere\thanks{D\'epartement de math\'ematiques, Universit\'e d'Angers, 2 Boulevard Lavoisier, 49045 Angers cedex 01 (France)}}

\date{}

\newtheorem{prop}[subsection]{Proposition}
\newtheorem{thm}[subsection]{Theorem}
\newtheorem{cor}[subsection]{Corollary}
\newtheorem{fact}[subsection]{Fact}
\newtheorem{conj}[subsection]{Conjecture}
\newenvironment{rem}%
   {\refstepcounter{subsection}%
        \medbreak\noindent{\bf Remark \thesubsection\space}}%
   {\par\medbreak}%
\newenvironment{exmp}%
   {\refstepcounter{subsection}%
        \medbreak\noindent{\bf Example \thesubsection\space}}%
   {\par\medbreak}%
\newenvironment{proof}%
   {\medbreak\noindent{\it Proof:\space}}%
   {\par\noindent\vrule height 5pt width 5pt depth 0pt\smallbreak}%

\newcommand{\df}{\bf}
\renewcommand{\cal}{\mathcal}
\let\sauvegardetiret=\-
\renewcommand{\-}[1]{\ifx#1-\penalty10000\hbox{-\relax}\penalty10000\else\sauvegardetiret#1\fi}
\newcommand{\Sp}{{\rm Spec}}
\newcommand{\A}{\forall}
\newcommand{\E}{\exists}
\newcommand{\tq}{\,\big/\ }
\newcommand{\nin}{\not\in}
\newcommand{\vect}{\overrightarrow}

\newcommand{\disj}{\mathop{\lor\mskip-6mu\relax\lor}}
\newcommand{\Lliff}{\;\longleftrightarrow\;}  
\newcommand{\Lssi}{\;\Longleftrightarrow\;}   
\newcommand{\donc}{\Rightarrow}          
\newcommand{\Ldonc}{\;\Longrightarrow\;} 
\newcommand{\Lto}{\;\longrightarrow\;}   
\newcommand{\NN}{{\mathbb N}}

\newcommand{\KK}{{K}}
\newcommand{\p}{{\rm\bf p}}
\newcommand{\ggf}{{\varphi}}
\newcommand{\cC}{{\cal C}}
\newcommand{\cI}{{\cal I}}
\newcommand{\cL}{{\cal L}}
\newcommand{\cN}{{\cal N}}
\newcommand{\cP}{{\cal P}}
\newcommand{\cU}{{\cal U}}
\newcommand{\cY}{{\cal Y}}
\newcommand{\UN}{{\rm\bf 1}}
\newcommand{\ZERO}{{\rm\bf 0}}
\newcommand{\meet}{\land}
\newcommand{\join}{\lor}
\newcommand{\mmeet}{\mathop{\meet\mskip-6mu\relax\meet}}
\newcommand{\jjoin}{\mathop{\join\mskip-6mu\relax\join}}
\newcommand{\cconj}{\mathop{\bigwedge\mskip-9mu\relax\bigwedge}}
\renewcommand{\disj}{\bigvee}
\newcommand{\ddisj}{\mathop{\bigvee\mskip-9mu\relax\bigvee}}
\newcommand{\lzar}{{\rm L_{Zar}}}
\newcommand{\llin}{{\rm L_{lin}}}
\newcommand{\ldef}{{\rm L_{def}}}

\newcommand{\lalin}{{\rm L_{lin}^{At}}}
\newcommand{\ladef}{{\rm L_{def}^{At}}}
\newcommand{\gen}[1]{\langle{#1}\rangle}
\newcommand{\szar}{{\rm SC_{Zar}}}
\newcommand{\slin}{{\rm SC_{lin}}}
\newcommand{\sdef}{{\rm SC_{def}}}
\newcommand{\sazar}{{\rm ASC_{Zar}}}
\newcommand{\salin}{{\rm ASC_{lin}}}
\newcommand{\sadef}{{\rm ASC_{def}}}
\newcommand{\AT}{{\rm At}}
\newcommand{\ASC}{${\rm ASC}$\-}
\newcommand{\asc}{\mathop{\rm asc}\nolimits}
\newcommand{\scdim}{\mathop{\rm sc\hbox{-}dim}\nolimits}
\newcommand{\pc}{{\rm C}}
\newcommand{\llat}{{\cal L}_{\rm lat}}
\newcommand{\ltc}{{\cal L}_{\rm TC}}
\newcommand{\lsc}{{\cal L}_{\rm SC}}
\newcommand{\lscd}[1]{{\cal L}_{{\rm SC}_{#1}}}
\newcommand{\lasc}{{\cal L}_{\rm ASC}}
\newcommand{\up}{\mathord\uparrow}
\begin{document}

\maketitle

\begin{abstract}
It is known from Grzegorczyk's paper \cite{grze-1951} that 
the lattice of real semi-algebraic closed subsets of ${\mathbb R}^n$ 
is undecidable for every integer $n\geq 2$. 
More generally, if $X$ is any definable set over a real 
or algebraically closed field $K$, then the lattice $L(X)$ 
of all definable subsets of $X$ closed in $X$ is undecidable whenever 
$\dim X\geq 2$. Nevertheless, we investigate in this paper 
the model theory of the class ${\rm SC_{def}}(K,d)$ of all such 
lattices $L(X)$ with $\dim X\leq d$ and $K$ as above or 
a henselian valued field of characteristic zero. 

We show that the universal theory of ${\rm SC_{def}}(K,d)$, 
in a natural expansion by definition of the lattice language, 
is the same for every such field $K$. We give a finite 
axiomatization of it and prove that it is locally finite 
and admits a model-completion, which turns to be decidable 
as well as all its completions. We expect $L({\mathbb Q}_p^d)$ 
to be a model of (a little variant of) this model-completion.
This leads us to a new conjecture in $p$-adic semi-algebraic 
geometry which, combined with the results of this paper, would 
give decidability (via a natural recursive axiomatization) 
and elimination of quantifiers for the complete theory 
of $L({\mathbb R}_p^d)$, uniformly in $p$. 
\end{abstract}

\section{Introduction}
\label{sec:introduction}

In this paper we study the model-theory of a class of lattices 
coming from the following examples. 

\begin{exmp}\label{exmp:Ldef}
Let $K$ be a henselian valued field of characteristic zero, 
a real closed field or an algebraically closed field. 
There exists a good notion of dimension for definables sets 
$A$ over $K$ (see \cite{Drie-1989} for the henselian case, 
and any book of real or complex algebraic geometry 
for the other cases). For any positive integer $i$ let: 
    $$\pc^i(A)=\overline{\{a\in A\tq \dim(A,a)=i\}}$$
where $\dim(A,a)$ is the maximal dimension of definable 
neighborhood of $a$ in $A$, and the overline stands for the 
topological closure in $A$. This is a definable subset of $A$
which we call the $i$\--pure\footnote{If $\pc^i(A)$ is non-empty 
it has pure dimension $i$, that is the local dimension $\dim(A,a)=i$
at every point of $A$.\label{fn:pure dim geom}}%
component of $A$. Given a definable set $X$ 
over $K$ of dimension $d$, let $\ldef(X)$ be the lattice of 
all definable subsets of $X$ closed in $X$, enriched\footnote{The 
additionnal functions are definable in the lattice structure 
of $\ldef(X)$.} with the unary functions $(\pc^i)_{i\leq d}$ 
and the binary function: 
$$A-B=\overline{A\setminus B}$$
Eventually let $\sdef(K,d)$ denote the class of lattices 
$\ldef(X)$ for all definable sets over $K$ of dimension 
at most $d$.  
\end{exmp}

\begin{exmp}
Let $K$ be any infinite field, $X$ a constructible subset 
over $K$ (that is a boolean combination of Zariski closed 
subsets of $K^n$ for some positive integer $n$) and 
$\lzar(X)$ be the lattice of all constructible subsets of 
$X$ which are Zariski closed in $X$ enriched with the 
following structure. For any $A,B$ in $\lzar(X)$ let $A-B$
be as in the above example and $C^i(A)$ be the union 
of the irreductible components of $A$ of Krull dimension $i$ 
(in the usual sense for topological spaces). 

Eventually let $\szar(K,d)$ denote the class of lattices 
$\lzar(X)$ for all definable sets over $K$ of dimension at most $d$.  
Of course $\szar(K,d)=\sdef(K,d)$ when $K$ is algebraically closed. 
\end{exmp}

It is known from an argument of \cite{grze-1951} that the complete 
theory of $\ldef(K^n)$ is undecidable for every real closed field $K$ 
and every integer $n\geq 2$, and the argument can easily be adapted 
to algebraically closed fields $K$. This paper gives some reason 
to believe that the complete theory of $\ldef(K^n)$ is decidable 
for every $p$-adically closed field $K$ and every $n$. It is 
organized as follows.

We first give in Section~\ref{sec:subscalled-lattices} 
a finite list of universal axioms of a theory $T_d$ in a language 
$\lscd d$ extending the language of lattices, the model of which 
we call $d$\--subscaled lattices. The examples given above 
are all $d$\--subscalled lattices. 
After some preliminar techinal results 
in Section~\ref{sec:embedd-subsc-latt} we prove 
in Section~\ref{sec:local-finiteness} that every finiteley 
generated $d$\--subscaled lattice is finite. 
Combining this result with a linear representation 
for finite $d$\--subscaled lattices and with the model-theoretic 
compactness theorem, we then prove in Section~\ref{sec:line-repr} 
that $T_d$ is precisely 
the universal theory of $\szar(\KK,d)$. In particular this 
theory is finitely axiomatizable and, remarkably enough, does 
not depend on $\KK$.
Eventually a detailed study of finitely generated 
extensions of finite $d$\--subscaled lattices, 
achieved in Sections~\ref{sec:primitive-extensions} 
and \ref{sec:signatures}
allows us to exhibit in Section~\ref{sec:model-completion} 
a model-completion $\bar T_d$ for $T_d$, 
having a finite axiomatisation. Moreover we show that 
$\bar T_d$ has finitely many completions, each of 
which is finitely axiomatizable and $\aleph_0$\--categorical. 
It follows that $\bar T_d$ is decidable.

It is difficult to find a model of $\bar T_d$ coming from 
geometry because such models are atomless. 
We present in the last section a similar model completion 
and decidability result for a theory $\bar T_d^*$ 
authorizing atoms. These results lead us to the following 
conjecture (or question): 

\begin{conj}
Let $K$ be a $p$\--adically closed field and $A$ be an infinite 
definable subset of $K^n$ which is open in its closure. 
Let $(B_k)_{k\leq q}$ be a finite collection of closed definable 
subsets of $\bar A\setminus A$. Then there exists a collection 
$(A_k)_{k\leq q}$ of non-empty definable subsets of $A$ clopen in $A$ 
such that: 
$$\A k\leq q,\quad \overline{A_k}=A_k\cup B_k$$
\end{conj}

If this conjecture is true then it follows immediatly that 
$\ldef(K^n)$ is a model of $T_n^*$ hence has a decidable 
complete theory (not depending on $p$).

\begin{rem}
Since $0$\--subscaled lattices are exactly non-trivial 
boolean algebras (with the $\lscd 0$\--structure and the boolean 
structure being quantifier-free bi-definable) our model-completion 
result for subscaled lattices is a generalisation to arbitrary finite 
dimension $d$ of the well known theorem on the model-completion 
of boolean algebras.
\end{rem}

\begin{rem}
The duals of $d$\--subscaled lattices form an elementary class 
of Heyting algebras so this paper may also be considered 
as a contribution to the model-theory of Heyting algebras. 
However the usual geometric objects whose study motivated 
this paper are closed sets (points, curves, surfaces, and so on). 
From this point of view the lattice $\lzar(\KK^d)$ is a more 
natural object to consider than its dual, the Heyting algebra 
of {\it open} algebraic sets. 
This is the reason why we had to present our results in this 
settings and not in terms of Heyting algebra.
\end{rem}

\section{Notation and definitions}
\label{sec:notation-definitions}
\label{sec:subscalled-lattices}

The set of all positive integers is denoted by $\NN$,
and $\NN^*=\NN\setminus\{0\}$. If $\cN$ is an unbounded subset
of $\NN$ (resp. the empty subset) we set $\max\cN=\infty$ 
(resp. $\max\cN=-1$). The symbols $\subseteq$ and $\subset$ 
denote respectively the inclusion and the strict inclusion.

\subsection{Lattices}

Let $\llat=\{\ZERO,\UN,\join,\meet\}$ be the language of 
lattices. An upper semi-lattice is an 
${\cal L}_{\rm up}$\--substructure of a lattice, 
with $ {\cal L}_{\rm up}=\{\ZERO,\join\}$. 
As usually $b\leq a$ is an abreviation for 
$a\join b =a$ and similarly for $b<a$, $b\geq a$ and 
$b>a$. Iterated $\join$ and $\meet$ operations are denoted 
by $\jjoin_{i\in I}a_i$ and $\mmeet_{i\in I}a_i$ respectively.
If the index set $I$ is empty then of course 
$ \jjoin_{i\in I}a_i=\ZERO$ and $\mmeet_{i\in I}a_i=\UN$.
The logical connectives `or', `and' and their iterated 
forms will be denoted by $\bigwedge$, $\disj$, $\cconj$ and $\ddisj$ 
respectively. We consider the following relation, definable in 
any lattice: 
\begin{eqnarray*}
b\ll a 
  &\iff& b<a \hbox{ and } \A c\;(c< a\donc b\join c < a) \\
  &\iff& b\leq a\neq\ZERO \hbox{ and } \A c\;(c< a\donc b\join c < a)
\end{eqnarray*}

The {\df spectrum} of a distributive lattice $L$ is 
the set $\Sp(L)$ of all prime filters of $L$, endowed 
with the so-called Zarisky topology, defined by taking
as a basis of closed sets the family:
\[
P(a)=\{\p\in\Sp(L)\tq a\in\p\}
\]
for $a$ ranging over $L$. Stone-Priestley's duality asserts
that $a\mapsto P(a)$ is an isomorphism between $L$ 
and the lattice of closed subsets of $\Sp(L)$ whose complement 
in $\Sp(L)$ is compact. We call a lattice {\df noetherian} 
if it is isomorphic to the lattice of closed sets of a noetherian 
topological space. By Stone-Priestley's duality a lattice $L$ 
is noetherian if and only if its spectrum is a noetherian 
topological space. In such a lattice every filter is principal 
and every element $a$ writes uniquely as the join of its 
{\df $\join$\--irreducible components}, which are the 
(finitely many) maximal elements in the set of 
non-zero $\join$\--irreducible elements of $L$ smaller than $a$. 
We denote by $\cI(L)$ the set of all non-zero 
$\join$\--irreducible elements of $L$.

We define the {\df lattice dimension
of an element $a$ in a lattice $L$} as the least upper 
bound (in $\NN\cup\{-1,\infty\}$) of the set of positive 
integers $n$ such that:
$$\E\p_0\subset\cdots\subset\p_n\in P(a)$$
This is nothing but the ordinary topological dimension 
(defined by chains of irreducible closed subsets) of 
the spectral space $P(a)$. 
We denote this dimension by $\dim_L a$. By construction 
$\dim_L a=-1$ if and only if $a=\ZERO_L$. The index $L$ 
is necessary since $\dim_L a$ is not preserved by 
$\llat$\--embeddings, nevertheless we omit it whenever
the ambiant lattice is clear from the context and 
we often call the lattice dimension simply the dimension. 
We let the {\df lattice dimension of $L$} 
be the lattice dimension of $\UN_L$ in $L$.

\begin{fact}\label{fact:ldim et dim geom}
If $L=L(X)$ is any of the lattices of the introduction, 
then for any $A\in L(X)$, $\dim_{L(X)} A$ is exactly the 
usual dimension of $A$ as a definable (or constructible) 
set over $K$. 
\end{fact}
Fact~\ref{fact:ldim et TCdim} below is a key argument 
in the proof of this result, which is non-obvious because 
the definition of $\dim_{L(X)} A$ lies on prime filters 
of closed definable subsets of $A$, an object which 
do not have a natural geometric meaning.

\subsection{TC-lattices}
Let $\ltc=\llat\cup\{-\}$ with `$-$' a binary 
function symbol. A {\df topologically complemented lattice},
or {\df TC-lattice} for short, is an $\ltc$\--structure 
which is a lattice and in which the {\df relative topological 
complement} $a-b$ is defined as the least element 
$c$ such that $a\leq b\join c$ (or equivalently $P(a-b)$ 
is the topolological closure of the relative complement 
$P(a)\setminus P(b)$, so the name). 
This is clearly the dual of a Heyting algebra with $a-b$ in the 
TC-lattice becoming $b\to a$ in its dual. So we know from 
the theory of Heyting algebra (see for example\cite{John-1982}) 
that every TC-lattice 
is distributive, and that the class of all TC-lattices 
is a variety (in the sense of universal algebra). 
Observe that in TC-lattices the $\ll$ relation is 
quantifier-free definable since: 
\[
a-b=a \iff b\meet a\ll a\hbox{ or }a=\ZERO
\]
So it will be preserved by $\ltc$\--embeddings. 
On the other hand the lattice dimension will not be 
preserved in general by $\ltc$\--embeddings 
of TC\--lattices.

We will use the following rules, the proof of which are
elementary exercises (using either dual properties, if known, 
of Heyting algebras or, more directly, Stone-Priestley's duality). 
\newcommand{\tcref}[1]{$\rm TC_{\ref{#1}}$}
\begin{list}{$\rm\bf TC_\theenumi:$}{\usecounter{enumi}}
\item\label{TC:a=(a inter b) union (a-b)}
$a=(a\meet b)\join(a-b)$.\\
In particular if $a$ is $\join$\--irreducible then 
$b< a\Ldonc b\ll a$.
\item\label{TC:(x union y)- z}
$(a_1\join a_2)-b = (a_1-b)\join (a_2-b)$.
\item\label{TC:(x-y)-y}
$(a-b)-b=a-b$.\\
So either $a-b=\ZERO$ or $(a-b)\meet b\ll a-b \leq a$.
\item\label{TC:x - (y union z)}
More generally $a-(b_1\join b_2)=(a-b_1)-b_2$.\\
So if $a-b_1=a$ then $a-(b_1\join b_2)=a-b_2$.
\end{list}

\begin{fact}\label{fact:ldim et TCdim}
For any TC\--lattice $L$ and any $a\in L$, 
$\dim_L a$ is exactly the least upper bound of the set 
of positive integers $n$ such that there exists
$a_0,\dots,a_n\in L$ such that: 
\[
\ZERO\neq a_0\ll a_1\ll\cdots\ll a_n\leq a
\]
\end{fact}
The proof is a good exercise that we leave to the reader. 
The fact that $\dim_L a$ is at least equal to the above 
least upper bound is true in any lattice. Equality
holds in TC\--lattices because they are min-compact (that 
is for any $a\in L$ the set of elements of $P(a)$ which 
are minimal with respect to the inclusion, is compact).

\subsection{(Sub)scaled lattices.}
For a given positive integer $d$ let 
$\lscd d=\ltc\cup\{\pc^i\}_{0\leq i\leq d}$ 
where the $\pc^i$'s are unary function symbols. 
With the examples of the introduction in mind, define 
the {\df sc\--dimension of an element $a$ of an 
$\lscd d$\--structure $L$} as: 
\[
\scdim a = 
\min\bigl\{l\leq d\tq a=\jjoin_{0\leq i\leq l} \pc^i(a)\bigr\}
\in\NN\cup\{-1,\infty\}
\]
The {\df sc\--dimension of $L$}, denoted $\scdim(L)$, is by 
definition the sc\--dimension of $\UN_L$.
In general the lattice dimension of an element is not 
preserved by $\llat$\--embeddings neither by $\ltc$\--lattices. 
On the other hand the sc\--dimension of 
an element is obviously preserved by $\lscd d$\--embeddings.

A {\df $d$\--subscaled lattice} is then an $\lscd d$\--structure  
which is a topologically complemented lattice and 
which satisfies the following list of axioms:
\newcommand{\scref}[1]{$\rm SC_{\ref{#1}}$}
\begin{list}{$\rm\bf SC_{\theenumi}:$}{\usecounter{enumi}}
\item
\label{SC:a=union des Ci(a)}
$\displaystyle
\jjoin_{0\leq i\leq d} \pc^i(a)=a$
\item
\label{SC:Ck(Ck(a))}%
\label{SC:Ck(Ci(a))}%
$\A I\subseteq\{0,\dots,d\},\ \A k$:
\[
\pc^k\Bigl(\jjoin_{i\in I} \pc^i(a)\Bigr)=
\left\{
\begin{array}{cl}
\ZERO  & \hbox{if } k\nin I\\
\pc^k(a) & \hbox{if } k\in I\\
\end{array}
\right.
\]
\item
\label{SC:Ck(Ck union Ck)}%
\label{SC:Ck(a union b)}%
$\A k\geq \max(\scdim(a),\scdim(b))$,\quad
$\pc^k(a\join b)=\pc^k(a)\join \pc^k(b)$
\item\label{SC:Ck(a inter b)}%
$\A i\neq j$,\quad 
$\scdim \left(\pc^i(a) \meet \pc^j(b) \right) < \min(i,j)$
\item
\label{SC:Ck(a) - b = Ck(a)}%
\label{SC:Ck(a) moins Cl(b)}%
\label{SC:Ck(Ck(a) moins Ck(b)}%
$\A k\geq \scdim(b)$,\quad $\pc^k(a)- b = \pc^k(a)- \pc^k(b)$\\
In particular, by \scref{SC:Ck(a union b)}: 
$\scdim b<a\Ldonc\pc^k(a)- b = \pc^k(a)$.
\item\label{SC:ll et scdim}
If $b\ll a$ then $\scdim b<\scdim a$.
\end{list}
\newcounter{sauvenumi}
\setcounter{sauvenumi}{\theenumi}

Axioms \scref{SC:a=union des Ci(a)} 
to \scref{SC:Ck(a) - b = Ck(a)} are easily seen to be 
equivalent to a finite set of equations in $\lscd d$. On the other 
hand \scref{SC:ll et scdim} is expressible by a universal 
formula in $\lscd d$ but not by an equation (indeed the 
class of $d$\--subscaled lattices is not preserved by 
$\lscd d$\--projections, hence is not a variety). 
We call {\df $d$\--scaled lattices}
the $d$\--subscaled lattices satisfying the following additional 
property: 
\begin{list}{$\rm\bf SC_{\theenumi}:$}{\usecounter{enumi}}
\setcounter{enumi}{-1}
\item
\label{SC:def scaled}
$\scdim a=\dim a$
\end{list}
This is an elementary class, which is not preserved 
by $\lscd d$\--substructures. By Fact~\ref{fact:ldim et dim geom} 
all the examples of $\lscd d$\--structures given in the introduction 
are $d$\--scaled lattices in which the lattice dimension and 
the sc\--dimension coincide with the usual geometric dimension.
As the terminology suggests, we will see that $d$\--subscaled 
lattices are precisely the $\lscd d$\--substructures of $d$\--scaled 
lattices. 
\medbreak

The finite language $\lscd d$ allows to write {\em finite} 
axiomatisations. When this is not essential we consider 
the language $\lsc=\ltc\cup\{\pc^i\}_{i\in\NN}$ and define 
{\df subscaled lattices} (resp. {\df scaled lattices}) 
as those $\lsc$\--structures whose $\lscd d$\--reduct is a 
$d$\--subscaled lattice (resp. a $d$\--scaled lattice)
for every $d$ large enough. This is not an elementary class, 
but for any fixed integer $d$ the class of (sub)scaled lattices 
of sc\--dimension at most $d$ (resp. exactly $d$) is elementary.  
Any $d$-(sub)scaled lattice expands uniquely to a (sub)scaled 
lattice of dimension at most $d$ by realizing $\pc^i$ as the 
constant map $x\mapsto\ZERO$ for every $i>d$. Conversely 
every (sub)scaled lattice $L$ is of that kind for every 
$d\geq\scdim(L)$.

\subsection{Basic properties}

The next additionnal properties follow easily
from the axioms of subscaled lattices.

\begin{list}{$\rm\bf SC_{\theenumi}:$}{\usecounter{enumi}}
\setcounter{enumi}{\thesauvenumi}
\item
\label{SC:scdim = max k}
$\scdim a=\max\{k\tq \pc^k(a)\neq\ZERO\}$\\
In particular $\A k,\ \scdim \pc^k(a)=k\Lssi\pc^k(a)\neq\ZERO$
\item
\label{SC:scdim a union b}%
$\scdim a\join b =\max(\scdim a, \scdim b)$\\
In particular $b\leq a\donc\scdim b\leq \scdim a$.
\item
\label{SC:determination de la main composante}%
$\displaystyle
\A k,\quad 
\pc^k(a)=\jjoin
  \Bigl\{ b\tq b\leq \jjoin_{0\leq i\leq k}\pc^i(a)
     \hbox{ and }\pc^k(b)=b \Bigr\}$
\item
\label{SC:scdim inferieure a dim}
$\scdim a\leq \dim a$
\item
\label{SC:a - Cd(a)}
$\displaystyle
\A I\subseteq\{0,\dots,d\},\quad
a-\jjoin_{i\in I}\pc^i(a)=\jjoin_{i\nin I}\pc^i(a)
$\\
In particular $\scdim (a-\pc^d(a))<d$.
\item
\label{SC:Ck(a) - b = Ck(Ck(a) - b)}
$\displaystyle
\A k,\quad 
\pc^k(a)-b=\pc^k(\pc^k(a)-b)
$
\item
\label{SC:Ck(a) est k pure}
$\displaystyle
\A k,\quad
\pc^k(a)=a \Lssi
\A b\ \big(a-b\neq\ZERO\donc\scdim a-b = k\bigr)
$
\end{list}

\begin{proof} (Sketch)
\scref{SC:scdim = max k} follows from 
\scref{SC:a=union des Ci(a)} and \scref{SC:Ck(Ck(a))};
\scref{SC:scdim a union b} from 
\scref{SC:Ck(Ck(a))}, \scref{SC:Ck(a union b)} 
and \scref{SC:scdim = max k}; 
\scref{SC:determination de la main composante} from
\scref{SC:Ck(Ck(a))} and \scref{SC:Ck(Ck union Ck)};
eventually \scref{SC:scdim inferieure a dim} is equivalent 
to \scref{SC:ll et scdim} modulo the other axioms.
Only the three last properties require a little effort. 
\smallskip

{\bf \scref{SC:a - Cd(a)}:} 
For every $l\in I$, $\pc^l(a)\leq \jjoin_{i\in I}\pc^i(a)$
hence $ \pc^l(a)-\jjoin_{i\in I}\pc^i(a)=\ZERO$. 
On the other hand for every $l\nin I$ and 
every $i\in I$, $\pc^l(a)-\pc^i(a)=\pc^l(a)$ 
by \scref{SC:Ck(a inter b)} and \scref{SC:Ck(Ck(a) moins Ck(b)}. 
So $\pc^l(a)-\jjoin_{i\in I}\pc^i(a)=\pc^l(a)$ 
by \tcref{TC:x - (y union z)}. Eventually by 
\scref{SC:a=union des Ci(a)} and \tcref{TC:(x union y)- z}:
\[
a-\jjoin_{i\in I}\pc^i(a)
=\jjoin_{l\leq d}\Bigl(\pc^l(a)-\jjoin_{i>k}\pc^i(a)\Bigr)
=\jjoin_{l\nin I}\pc^l(a)
\]
In particular if $I=\{d\}$ then 
$a-\pc^d(a)=\jjoin\limits_{k<d}\pc^k(a)$ 
hence $\scdim(a-\pc^d(a))<d$.
\smallskip

{\bf \scref{SC:Ck(a) - b = Ck(Ck(a) - b)}:}
By \tcref{TC:a=(a inter b) union (a-b)}, 
$\pc^k(a)=(\pc^k(a)-b)\join(\pc^k(a)\meet b)$, and by 
\scref{SC:scdim a union b} both $\pc^k(a)-b$ and $\pc^k(a)\meet b$ 
have sc\--dimension at most $k$, so by \scref{SC:Ck(Ck union Ck)}:
\[
\pc^k(\pc^k(a))=\pc^k(\pc^k(a)-b)\join\pc^k(\pc^k(a)\meet b)
\]
By \scref{SC:a=union des Ci(a)}, 
$\pc^k(\pc^k(a)\meet b)\leq\pc^k(a)\meet b\leq b$, 
and by \scref{SC:Ck(Ck(a))},  $\pc^k(a)=\pc^k(\pc^k(a))$ so:
\[
\pc^k(a)\leq \pc^k(\pc^k(a)-b)\join b\]
It follows that $\pc^k(a)-b\leq \pc^k(\pc^k(a)-b)$, 
and equality holds by \scref{SC:a=union des Ci(a)}.

\smallskip
{\bf \scref{SC:Ck(a) est k pure}:}
If $a=\pc^k(a)$ then $a-b=\pc^k(a-b)$ for every $b$ 
by \scref{SC:Ck(a) - b = Ck(Ck(a) - b)}. If moreover 
$a-b\neq\ZERO$ then $\scdim(a-b)=k$ by \scref{SC:scdim = max k}.
Conversely assume that $a\neq\pc^k(a)$ (hence $a\neq\ZERO$). 
If $\scdim a\neq k$ then $a-\ZERO=a\neq\ZERO$ hence
$\scdim a-\ZERO\neq k$. If $\scdim a=k$, let $b=\pc^k(a)$. 
Then by assumption $a-b\neq \ZERO$, and by \scref{SC:a - Cd(a)}, 
$a-b=\jjoin_{i<k} \pc^i(a)$ hence $\scdim a-b<k$.
\end{proof}

\subsection{Miscellanies}

Given an integer $k$ we say that an element $a$ 
of a distributive lattice $L$ is {\df $k$\--pure in $L$} 
if and only if\footnote{Remember Example~\ref{exmp:Ldef} and 
Footnote~\ref{fn:pure dim geom}.}:
\[
\A b\in L\ (a-b\neq\ZERO\donc \dim_L a-b =k)
\]
Then either $a=\ZERO$ or $\dim_L a = k$. In the latter case 
we say that $a$ {\df has pure dimension $k$ in $L$}. 

Similarly in any $d$\--subscaled lattice we say that an 
element $a$ is {\df $k$\--sc\--pure} if and only if:
\[
\A b\in L\ (a-b\neq\ZERO\donc \scdim a-b =k)
\]
By \scref{SC:Ck(a) est k pure}, $a$ is $k$\--sc\--pure 
if and only if $a=\pc^k(a)$. Then by~\scref{SC:scdim = max k}, 
either $a=\ZERO$ or $\scdim a=k$. In the latter case we 
say that $a$ {\df has pure sc\--dimension $k$}. 
For any $a$, the element $\pc^k(a)$ is called the 
{\df $k$\--sc\--pure component of $a$}, 
and simply its {\df $k$\--pure component} 
if $L$ is a scaled lattice.

\begin{prop}\label{prop:scaled structure definissable}
The $\lscd d$\--structure of a $d$\--scaled lattice $L$ is 
uniformly definable in the $\llat$\--structure of $L$. 
In particular it is uniquely determined by this 
$\llat$\--structure. 
\end{prop}

\begin{proof}
Clearly the TC-structure is an extension by definition 
of the lattice structure of $L$. For every positive 
integer $k$ the class of $k$-pure elements is uniformly
definable, using the definablility of $\ll$. Then  
so is the function $\pc^k$ for every $k$, by decreasing 
induction on $k$. 
Indeed by \scref{SC:determination de la main composante}
and \scref{SC:a - Cd(a)}, $\pc^k(a)$ is the largest 
$k$\--pure element $c$ such that $c\leq a-\jjoin_{i>k}\pc^i(a)$.
\end{proof}

Our study of (sub)scaled lattices is motivated by the examples 
given in the introduction. Although they are less natural, 
the following examples which will be needed further in this paper.

\begin{exmp}\label{exple:Lnoeth}
In an arbitrary noetherian lattice $L$ an $\lsc$\--structure 
can be defined as follows. For every $a,b\in L$, if $\cC(a)$ denotes 
the set of all $\join$\--irreducible components of $a$, let:
\begin{eqnarray*}
&a-b=\jjoin\{c\in\cC(a)\tq c\not\leq b\}& \\
&(\A k)\quad \pc^k(a)=\jjoin\{c\in\cC(a)\tq\dim_L c=k\}&
\end{eqnarray*}
This $\ltc$\--structure (resp. $\lsc$\--structure) is the only one
(by Proposition~\ref{prop:scaled structure definissable})
which turns $L$ into a TC\--lattice (resp. a scaled lattice). 
On the other hand, for any strictly increasing map 
$D\colon\cI(L)\to\NN$ and any $a,b\in L$ define $a-b$ as above
and: 
\[
(\A k)\quad \pc_D^k(a)=\jjoin\{c\in\cC(a)\tq D(c)=k\}
\]
This $\lsc$\--structure turns $L$ into a subscaled lattice 
which is not a scaled lattice (except if $D$ coincides with 
the map $\dim_L$). We will use without further mention 
the following obvious fact: 
\end{exmp}

\begin{fact}
Every noetherian (hence in particular
every finite) subscaled lattice is of the above kind.
\end{fact}

Eventually the following notation will be convenient in 
induction arguments. If $\cL$ is any of our languages 
$\llat$, $\ltc$, $\lscd d$ or $\lsc$ we let 
$\cL^*=\cL\setminus\{\UN\}$. Given an $\cL$\--structure
$L$ whose reduct to $\llat$ is a lattice, for any $a\in L$
we denote by: 
\[L(a)=\{b\in L\tq b\leq a\}\]
$L(a)$ is a typical example of $\cL^*$\--substructure of $L$.

\section{Embeddinds of subscaled lattices}
\label{sec:embedd-subsc-latt}

We need a reasonably easy criterion for an 
$\llat$\--embedding of subscaled lattices to be an 
$\lsc$\--embedding. 
The special case of a noetherian embedded lattice, presented 
in the next proposition, is sufficient for this paper. 
However combining the model-theoretic compactness theorem 
with the local finiteness Theorem~\ref{thm:TCS-treillis type fini}, 
one can easily derive from 
Proposition~\ref{prop:CNS de lsc plongement} that 
an $\llat$\--embedding $\varphi:L\to L'$ between arbitrary 
subscaled lattices is an $\lsc$\--embedding if and only if 
it preserves the sc\--dimension and sc\--purity, that is for 
every $a\in L$ and every $k\in\NN$:
\[
\pc^k(a)=a\Ldonc \pc^k(\varphi(a))=\varphi(a)
\]

\begin{prop}%
\label{prop:CNS de lsc plongement}%
Let $L_0$ be a noetherian subscaled lattice and $\varphi:L_0\to L$ 
an $\llat$\--embedding such that for every $a\in\cI(L_0)$, 
$\varphi(a)$ is sc\--pure and has the same sc\--dimension as $a$. 
Then $\varphi$ is an $\lsc$\--embedding.
\end{prop}

\begin{rem}
Clearly the same statement remains true with $\llat$ and $\lsc$ 
replaced respectively by $\llat^*$ and $\lsc^*$ (or $\lscd d^*$). 
We will freely use these variants.
\end{rem}

\begin{proof}
Let $d=\scdim L$, and for any $a\in L_0$: 
\[
(a_0,\dots,a_d)=(\varphi(\pc^0(a)),\dots,\varphi(\pc^d(a)))
\]
For every positive integer $k$, $\pc^k(a)$ 
is $k$\--sc\--pure by \scref{SC:Ck(a) est k pure}, hence each 
$\join$\--irreducible component $c$ of $\pc^k(a)$ in $L_0$ 
has pure sc\--dimension $k$. 
By assumption each such $\varphi(c)$ then has pure 
sc\--dimension $k$. 
The join of finitely many elements of pure sc\--dimension 
$k$ is easily seen to be $k$\--sc\--pure by definition 
and by \tcref{TC:(x union y)- z}, so we have proved:
\begin{equation}
\A k,\ a_k\hbox{ is $k$-sc-pure.} \label{eq:phi k pure}
\end{equation}

Moreover for any $k\neq l$, $\scdim(\pc^k(a)\meet\pc^l(a))<\min(k,l)$ 
by \scref{SC:Ck(a inter b)}. It follows that each 
$\join$\--irreducible component $c$ of $\pc^k(a)\meet\pc^l(a)$ 
has sc\--dimension stricly less than $\min(k,l)$, hence so does 
$\varphi(c)$ by assumption. By \scref{SC:scdim a union b} 
we conclude that:
\begin{equation}
\A k\neq l,\ \scdim(a_k\meet a_l)<\min(k,l) 
\label{eq:dim de phi ak inter al}
\end{equation}

For every $k>d$, $\varphi(\pc^k(a))=\varphi(\ZERO)=\ZERO$. 
Each $\join$\--irreducible component of $a$ has sc\--dimension 
at most $d$ hence so does $\varphi(a)$ by assumption, so by 
\scref{SC:scdim a union b}, $\scdim(\varphi(a))\leq d$. 
It follows that $\pc^k(\varphi(a))=\ZERO=\pc^k(\varphi(a))$. 
\smallskip

For every $k\leq d$ let: 
\[
b_k=\jjoin_{0\leq l\leq k}a_l
\qquad\hbox{and}\qquad
c_k=\jjoin_{0\leq l\leq k}\pc^l(\varphi(a))
\]
By \scref{SC:a - Cd(a)}, $c_k-\pc^k(\varphi(a))=c_{k-1}$. 
Moreover the proof of \scref{SC:a - Cd(a)} proves as well that 
$b_k-a_k=\jjoin_{l\leq k}(a_l-a_k)=b_{k-1}$, thanks 
to (\ref{eq:phi k pure}) and (\ref{eq:dim de phi ak inter al})
above. 
\smallbreak

Now assume that $b_k=c_k$ for some $k\leq d$. 
Then $\pc^k(\varphi(a))=\pc^k(c_k)$ by \scref{SC:Ck(Ck(a))}, 
so $\pc^k(\varphi(a))=\pc^k(b_k)$ by assumption. 
It follows by \scref{SC:Ck(Ck(a) moins Ck(b)} that:
\[
\pc^k(\varphi(a))-a_k =
\pc^k(b_k)-\pc^k(a_k) =
\pc^k(b_k-a_k) =
\pc^k(b_{k-1})
\]
\scref{SC:scdim a union b} implies that $\scdim(b_{k-1})\leq k-1$ 
so $\pc^k(b_{k-1})=\ZERO$ by \scref{SC:Ck(Ci(a))}.
It follows that $\pc^k(\varphi(a))-a_k=\ZERO$ that is 
$\pc^k(\varphi(a))\leq a_k$. 
On the other hand $a_k\leq \jjoin_{l\leq k} \pc^k(\varphi(a))$ 
and $\pc^k(a_k)=a_k$ by (\ref{eq:phi k pure}) so 
$a_k\leq \pc^k(\varphi(a))$ by 
\scref{SC:determination de la main composante}.
We have proved: 
\[
b_k=c_k\Ldonc a_k= \pc^k(\varphi(a))
\]
Then $b_{k-1}=b_k-a_k=c_k-\pc^k(\varphi(a))=c_{k-1}$. 
Since $b_d=c_d=a$ it follows by decreasing induction that 
$b_k=c_k$ for every $k$, that is 
$\varphi(\pc^k(a))=\pc^k(\varphi(a))$.
Incidentally, since $\varphi$ is injective, this implies 
by \scref{SC:scdim = max k} that for every $a\in L_0$: 
\begin{equation}\label{eq:phi preserve scdim}
\scdim a=\scdim\varphi(a)
\end{equation}

Now let $a,b\in L_0$, and $a',b'$ be their images by $\varphi$. 
We have to show that $\varphi(a-b)=a'-b'$. 
By \tcref{TC:(x union y)- z}, replacing if necessary $a$ 
by its $\join$\--irreducible components, we may assume w.l.o.g. 
that $a$ itself is $\join$\--irreducible in $L_0$. This implies 
that $a=\pc^k(a)$ for some $k$. 
It then remains two possibilities for $a-b$:
\begin{itemize}
\item 
If $b\geq a$ then $a-b=\ZERO$ and $b'\geq a'$ hence
$\ggf(a-b)=\ZERO=a'-b'$.
\item
Otherwise $c=b\meet a<a$ hence $c\ll a$ by 
\tcref{TC:a=(a inter b) union (a-b)}, so $\scdim c<\scdim a$
by \scref{SC:ll et scdim}. 
Let $c'=\ggf(c)=b'\meet a'$, by assumption 
$a'=\pc^k(a)$ and by (\ref{eq:phi preserve scdim}) 
$\scdim(c')<\scdim(a')=k$,
hence by \scref{SC:Ck(Ck(a) moins Ck(b)} $a'-c'=a'$. 
We conclude that $a-b=a$ and $a'-b'=a'$ hence
$\ggf(a-b)=\ggf(a)=a'=a'-b'$.
\end{itemize}
We have proved that $\varphi$ preserves $-$ and the 
$\pc^k$'s operations, so $\varphi$ is an $\lsc$\--embedding.
\end{proof}

\begin{cor}\label{cor:CNS ss treillis ss LSC structure}
Let $L_0$ be a noetherian sublattice of a subscaled lattice $L$. 
If for any $b<a\in\cI(L)$, $a$ is sc\--pure in $L$ and 
$\scdim b<\scdim a$ in $L$, then $L_0$ is an $\lsc$\--substructure 
of $L_0$. 
\end{cor}

\begin{proof}
The assumptions imply that the map $D:a\mapsto\scdim a$ 
is a stricly inreasing map from $\cI(L_0)$ to $\NN$. Endow 
$L_0$ with the structure of subscaled lattice determined 
by $D$ as in Example~\ref{exple:Lnoeth}. 
Proposition~\ref{prop:CNS de lsc plongement} then applies 
to the inclusion map $\varphi$ from $L_0$ to $L$. 
\end{proof}

\section{Local finiteness}
\label{sec:local-finiteness}

We prove in this section that every finitely generated subscaled 
lattice is finite. This result is far non-obvious, due to the lack 
of any known normal form for terms in $\lsc$. It contrasts with 
the situation in TC-lattices, which can be both infinite and 
generated by a single element. Our main ingredient, which explains 
this difference, is the uniform bound given {\it a priori} for the 
sc\--dimension of any element in a given subscaled lattice.

\begin{thm}\label{thm:TCS-treillis type fini}
Any subscaled lattice $L$ of sc\--dimension $d$ generated by $n$ 
elements is finite. More precisely, the cardinality of $\cI(L)$ 
is then bounded by the fonction $\mu(n,d)$ defined by: 
\[
\begin{array}{rcll}
\mu(n,-1) &=& 0    &\ (\A n)\\
\mu(n,d)  &=& 2^n+ \mu(2^{n+1},d-1) &\ (\A n,\ \A d\geq 0)
\end{array}
\]
\end{thm}

\begin{proof}
If $d=-1$ the only subscaled lattice of dimension $-1$ is 
the one-element lattice $\{\ZERO\}$, so the result is trivial. 

Assume the $d\geq 0$ and that the result is proved for every 
$d'<d$ and every positive integer $n$. Let $L$ be a subscaled 
lattice of sc\--dimension $d$ generated by elements 
$x_1,\dots,x_n$.
Let $\Omega_n$ be the 
family of all subsets of $\{1,\dots,n\}$ 
(so $\Omega_0=\{\emptyset\}$).  
For every $I\in \Omega_n$ 
let $I^c=\Omega_n\setminus I$ and:
\[
y_I = \Bigl(\mmeet_{i\in I} x_i\Bigr)
      - \Bigl(\jjoin_{i\in I^c} x_i\Bigr),
\qquad
z_I =\pc^d(y_I)
\]
The family of all 
$\cY_I = \bigcap_{i\in I} P(x_i) \cap \bigcap_{i\in I^c} P(x_i)^c$
is a partition of $\Sp(L)$. Indeed the $\cY_i$'s are the atoms 
of the boolean algebra generated in the power set $\cP(\Sp(L))$ 
by the $P(x_i)$'s. Moreover each $P(y_I)$ is the 
topological closure $\overline{\cY}_I$ of $\cY_I$ in $\Sp(L)$ 
hence for every $x\in L$:
\[
P(x)=\bigcup_{I\in\Omega_n}P(x)\cap \cY_I
\subseteq\bigcup_{I\in\Omega_n}P(x)\cap \overline{\cY}_I
=P\Bigl(\jjoin_{I\in\Omega_n}x\meet y_I\Bigr)
\]
So $x\leq \jjoin\limits_{I\in\Omega_n}(x\meet y_I)$ 
by Stone-Priestley's duality. The reverse inequality being 
obvious we have proved:
\begin{equation}
\label{eq:x et les xi}
\A x\in L,\quad x\leq \jjoin\limits_{I\in\Omega_n}(x\meet y_I)
\end{equation}
In particular \scref{SC:Ck(Ck union Ck)} also gives: 
\begin{equation}
\label{eq:Cd 1 egale union des Cd zI}
\pc^d(\UN)
=\pc^d\Bigl(\jjoin\limits_{I\in\Omega_n}y_I\Bigr)
=\jjoin\limits_{I\in\Omega_n}z_I
\end{equation} 
For every $I\neq J\in\Omega_n$, if for example $I\not\subseteq J$ 
choose any $i\in I\setminus J$ and observe that 
$y_I\leq x_i$ and $y_J\leq \UN-x_i$ so $y_I\meet y_J\ll\UN-x_i$
by \tcref{TC:(x-y)-y}. By \scref{SC:ll et scdim} and the 
$d$\--sc\--purity of the $z_I$'s it follows that:
\begin{equation}
\label{eq:zI meet zJ}
\scdim z_I\meet z_J<d\hbox{\quad hence\quad}z_I-z_J=z_I
\end{equation}
It follows from \scref{SC:scdim a union b}, \scref{SC:a - Cd(a)} 
and (\ref{eq:zI meet zJ}) above, that the element:
\[
a=\big(\UN-\pc^d(\UN)\big)\join
\Big(\jjoin\limits_{I\in\Omega_n}(y_I-z_I)\Big)\join
\Big(\jjoin\limits_{I\neq J\in\Omega_n}(z_I\meet z_J)\Big)
\] 
has sc\--dimension strictly smaller than $d$. So the 
induction hypothesis applies to the $\lsc$\--substructure 
$L_0^-$ of $L(a)$ 
generated by the ${(y_I-z_I)}$'s and the ${(z_I\meet a)}$'s: 
$L_0^-$ is finite, with at most $\mu(2|\Omega_n|,d-1)$ 
non-zero $\join$\--irreducible elements. 
Eventually let $L_1$ be the upper semi-lattice generated 
in $L$ by $L_0^-\cup\{z_I\}_{I\in\Omega_n}$. 
By construction $L_1$ is finite and 
$\cI(L_1)\subseteq\cI(L_0)\cup\{z_I\}_{I\in\Omega_n}$, 
so $|\cI(L_1)|\leq 2^n+\mu(2^{n+1},d-1)=\mu(n,d)$.
It is then sufficient to show that $L_1=L$. 
\smallskip

By~(\ref{eq:Cd 1 egale union des Cd zI}), 
$\UN=\pc^d(\UN)\join a=\jjoin_{I\in\Omega_n}z_I\join a\in L_1$.
For every $I\in\Omega_n$ and every $b\in L_0^-$, 
$z_I\meet b=(z_I\meet a)\meet b\in L_0^-$. 
Eventually for every $I\neq J\in\Omega_n$, 
$z_I\meet z_j=(z_I\meet a)\meet(z_J\meet a)\in L_0^-$. 
So by the distributivity law, $L_1$ is a sublattice of $L$. 

Since $\cI(L_1)\subseteq\cI(L_0^-)\cup\{z_I\}_{I\in\Omega_n}$ 
it is immediate that for any $b'<b$ in $\cI(L_1)$, 
$\scdim b'<\scdim b$. 
So $L_1$ is an $\lsc$\--sub\-structure of $L$ by 
Corollary~\ref{cor:CNS ss treillis ss LSC structure}.
Moreover each $y_I=(y_I-z_I)\join z_I\in L_1$ and for every $i\leq n$, 
(\ref{eq:x et les xi}) gives:
\[
x_i=\jjoin_{I\in\Omega_n}x_i\meet y_I
\leq \jjoin_{{I\in\Omega_n}\atop{i\in I}}y_I\leq x_i
\]
So equality holds, hence each $x_i\in L_1$ and eventually 
$L=L_1$. 
\end{proof}

\begin{cor}
\label{cor:structures a n generateurs}
For every $n,d$ there are finitely many non-isomorphic 
subscaled lattices of sc\--dimension $d$ generated by $n$ elements.
\end{cor}

\begin{proof}
Any such subscaled lattice $L$ is finite, with 
$|\cI(L)|\leq \mu(n,d)$ by Theorem~\ref{thm:TCS-treillis type fini}. 
Clearly there are finitely many non-isomorphic lattices such that 
$|\cI(L)|\leq \mu(n,d)$ and each of them admits finitely 
many non-isomorphic $\lscd d$\--structures. 
The conclusion follows.
\end{proof}

\section{Linear representation}
\label{sec:line-repr}

In this section we prove that the theory of $d$\--subscaled 
lattices is the universal theory of various natural classes 
of $\lscd d$\--structures, including $\szar(\KK,d)$. 
The argument is based on an elementary representation 
theorem for $d$\--subscaled lattices, combined with the 
local finiteness result of Section~\ref{sec:local-finiteness}.

Given an arbitrary field $\KK$, a non-empty linear variety 
$X\subseteq \KK^m$ is determined 
by the data of an arbitrary point $P\in X$ and 
the vector subspace $\vect X$ of $\KK^m$, 
{\it via} the relation $X=P+\vect X$ (the orbit of $P$ 
under the action of $\vect X$ by translation). 
We call $X$ a {\df special linear variety} (resp. a 
{\df special linear set}) if $X$ is a linear 
variety such that $\vect X$ is generated 
by a subset of the canonical basis of $\KK^m$ 
(resp. if $X$ is a finite union of special linear varieties). 
For example the empty set is a special linear set, as 
the union of an empty family of special linear varieties.
The family $\llin(X)$ of all special linear subsets of $X$ 
is a noetherian lattice (because it is the family of closed 
sets of a noetherian topology on $X$). So it has a natural 
structure of scaled lattice defined as in Example~\ref{exple:Lnoeth}. 

\begin{rem}\label{rem:sclin inclus dans scdef et sczar}
For every $A\in\llin(X)$, $\scdim A=\dim_{\llin(X)} A=$ 
the dimension of $A$ as defined in linear algebra. If $\KK$ is 
infinite then this dimension coincides with the Krull 
dimension of $A$. Moreover if $A$ is $\join$\--irreducible 
in $\llin(X)$ then it is pure dimensionnal, hence it is 
sc\--pure both in $\llin(X)$ and $\lzar(X)$. By
Proposition~\ref{prop:CNS de lsc plongement} it follows 
that if $\KK$ is infinite then $\llin(X)$ is an 
$\lsc$\--substructure of $\lzar(X)$. Similarly if $\KK$ is 
a henselian valued field of characteristic zero, a real closed 
field or an algebraically closed field then $\llin(X)$ is an 
$\lsc$\--substructure of $\ldef(X)$. 
\end{rem}

In the following proposition $\KK^m$ is identified 
to the subset $\KK^m\times\{0\}^{r}$ of $\KK^{m+r}$.

\begin{prop}%
\label{prop: treillis des lineaires speciaux}%
Given any two special linear sets 
$C\subseteq B\subseteq \KK^m$ and any $N\geq\dim C$ 
there exists a positive integer $r$ and a special 
linear set $A\subseteq \KK^{m+r}$ 
of pure dimension $n$ such that $A\cap B=C$. 
\end{prop}

\begin{proof}
For any integer $n$ let $(e_1,\dots,e_n)$ be the canonical 
basis of $\KK^n$. If $I$ is a subset of $\{1,\dots,n\}$ we 
denote $\vect E(I)$ the vector subspace of $\KK^n$ 
generated by $(e_i)_{i\in I}$. 

Decompose $C$ as a union of special linear varieties:
$C_1,\dots,C_p$ with each $C_i=P_i+\vect E(J_i)$ 
and $|J_i|=\dim C_i\leq n$. Let $r$ be larger than $0$, $N-m$ and 
every $N-|J_i|$. For each $i$ take an arbitrary 
subset $J'_i$ of cardinality $N-|J_i|$ inside 
$\{m+1,\dots, m+r\}$ and let  $I_i=J_i\cup J'_i$.
Choose an arbitrary point $P_0\in\KK^{m+r}\setminus \KK^m$
and let $I_0=\{1,\dots,N\}$. Eventually let 
$A_i=P_i+\vect E(I_i)$ for every $i\leq p$, 
and $A=A_0\cup\dots\cup A_p$. By construction each $A_i$ 
has dimension $|I_i|=N$, hence $A$ is $N$\--pure. 
At least $A_0$ is non empty (it is the only one 
in case $p=0$ that is if $C$ is empty) hence $A$ has pure 
dimension $N$. Clearly $A\cap\KK^m=C$ hence {\it a fortiori}
$A\cap B=C$.
\end{proof}

\begin{prop}[Linear representation]\label{prop:linear repres}
Let $\KK$ be an infinite field and let $L$ be any finite 
subscaled lattice. 
Then there exists a linear set $X$ over $\KK$ 
and an $\lsc$\--embedding $\varphi\colon L\to\llin(X)$. 
\end{prop}

\begin{rem}
Since $\varphi$, in the above proposition, is an $\lsc$\--embedding 
it preserves the sc\--dimension hence $\dim X=\scdim L$. 
So if $\scdim L\leq d$ we can indentify $L$ and $X$ with their 
respective $\lscd d$\--reduct, and $\varphi$ is then 
also an $\lscd d$\--embedding.
\end{rem}

\begin{proof}
Remember that $\lsc^*=\lsc\setminus\{\UN\}$. 
We prove by induction on the number $n$ of non-zero 
$\join$\--irreducible elements of a subscaled lattice $L$ 
that there exists an $\lsc^*$\--embedding $\varphi$ of $L$ 
into $\llin(\KK^m)$ for some $m$ depending on $L$. 
\smallbreak 

For $n=0$, $L$ is the one-element lattice $\{\ZERO\}$ 
hence it is an $\lsc^*$\--substructure of $\llin(K)$.

Let $n\geq 1$, assume the result proved for $n-1$
and take a subscaled lattice $L$ with non-zero 
$\join$\--irreducible elements $a_1,\dots,a_n$. 
Reordering if necessary we may assume that $a_n$ is maximal 
among the $a_i$'s. Let $a=a_n$, $b=\jjoin_{1\leq i<n}a_i$, 
$c=a\meet b$ and $\varphi$ an $\lsc^*$\--embedding of $L(b)$ 
into some $\llin(\KK^m)$ given by induction hypothesis.
Since $a$ is $\join$-irreducible in $L$ it is sc\--pure.
Moreover $c\ll a$ by \tcref{TC:a=(a inter b) union (a-b)}, 
hence $a$ has pure sc\--dimension $N$ for some $N>\scdim(c)$
by \scref{SC:ll et scdim}. 
Let $B,C$ be the respective images of $b,c$ by $\varphi$. 
Proposition~\ref{prop: treillis des lineaires speciaux}
gives a positive integer $r$ and a special linear set 
$A\subseteq \KK^{m+r}$ of pure dimension $N$ such that 
$A\cap B=C$. Since we indentified $\KK^m$ with 
$\KK^m\times\{0\}^r$ we consider $\llin(\KK^m)$ 
as an $\lsc^*$\--substructure of $\llin(\KK^{m+r})$. 
So $\varphi$ actually embeds $L(b)$ into $\llin(\KK^{m+r})$.

Every element $x$ of $L$ writes uniquely $x_a\join x_b$ 
with $x_a\in\{\ZERO,a\}$ and $x_b\in L(b)$ (group appropriately 
the $\join$\--irreducible components of $x$, using the maximality 
of $a_n$) hence:
\[
\bar\varphi(x)=
\left\{
\begin{array}{cl}
\varphi(x_b)      & \hbox{if $x_a=\ZERO$}\\
A\cup\varphi(x_b) & \hbox{if $x_a=a$}
\end{array}
\right.
\]
is a well-defined $\llat^*$\--embedding of $L$ into 
$\llin(\KK^{m+r})$. 
Moreover $\bar\varphi$  is an $\lsc^*$\--embedding by 
Proposition~\ref{prop:CNS de lsc plongement}.
This finishes the induction. 

We have constructed an $\lsc^*$\--embedding $\varphi$ of $L$ 
into $\llin(\KK^m)$ for some $m$. Then $X=\varphi(\UN_L)$ 
is a special linear set, so $\varphi$ induces 
an $\lsc$\--embedding of $L$ into $\llin(X)$.
\end{proof}

For any infinite field $\KK$ and positive integer $d$ 
let $\slin(\KK,d)$ be the class of $d$-scaled lattices $\llin(X)$ 
for every special linear variety $X$ over $\KK$ of dimension 
at most $d$. 

\begin{thm}
For any infinite field $\KK$ and positive integer $d$, 
the universal theories of $\slin(\KK,d)$ and $\szar(\KK,d)$ 
are exactly the theory of $d$-subscaled lattices. 

If moreover $\KK$ is a henselian valued field of characteristic 
zero, a real closed field or an algebraically closed field then 
the same holds for $\sdef(\KK,d)$.
\end{thm}

\begin{proof}
By Remark~\ref{rem:sclin inclus dans scdef et sczar} it 
suffices to prove the theorem for $\slin(\KK,d)$. Let 
$T(\KK,d)$ be its universal theory.

The linear representation Proposition~\ref{prop:linear repres} 
shows that every finite $d$\--subscaled lattice embeds into 
some $L\in\slin(\KK,d)$ hence is a model of $T(\KK,d)$. 
Since every finitely generated $d$\--subscaled lattice 
is finite by Theorem~\ref{thm:TCS-treillis type fini}, 
the model-theoretic compactness argument then implies 
that any $d$\--subscaled lattice is a model of $T(\KK,d)$.

Conversely every $L\in\slin(K,d)$ is a $d$\--scaled 
lattice hence obviously its universal theory contains 
the theory of $d$\--subscaled lattices.
\end{proof}

\section{Primitive extensions}
\label{sec:primitive-extensions}

This section is devoted to the study of minimal
proper extensions of finite subscaled lattices. 
Let $L_0$ be a finite subscaled lattice, $L$ an 
$\lsc$\--extension of $L_0$ and $x\in L$. 
We introduce the following notation.
\begin{itemize}
\item
For every $a\in L_0$, $a^-=\jjoin\{b\in L_0\tq b<a\}$.
\item
$L_0\gen x$ denotes the $\lsc$\--substructure of $L$ 
generated by $L_0\cup\{x\}$.
\item
$g(x,L_0)=\mmeet\{a\in L_0\tq x\leq a\}$. 
\end{itemize}
Clearly $a\in\cI(L_0)$ if and only if $a^-$ is the unique 
predecessor of $a$ in $L_0$ (otherwise $a^-=a$). 
We say that a tuple $(x_1,x_2)$ of 
elements of $L$ is {\df primitive over $L_0$} if there 
exists $g\in\cI(L_0)$ such that: 
\begin{enumerate}
\item 
$x_1$, $x_2$ are sc\--pure of the same sc\--dimension.
\item 
Each $g^-\meet x_i\in L_0$.
\item
One of the following happens: 
\begin{itemize}
\item
$x_1=x_2$ and $g^-\meet x_1\ll x_1 \ll g$.
\item
$x_1\neq x_2$, $x_1\meet x_2\in L_0$ and $g-x_1=x_2$, $g-x_2=x_1$.
\end{itemize}
\end{enumerate}
The above conditions imply that each $x_i\nin L_0$ and: 
\[
g=g(x_1,L_0)=g(x_2,L_0)
\]
We say that $L$ is {\df primitively generated over $L_0$}, 
or simply that it is a {\df primitive extension} of $L_0$,
if there exists $(x_1,x_2)$ primitive over $L_0$ such that 
$L=L_0\gen{x_1,x_2}$ (then clearly $L=L_0\gen x_1=L_0\gen x_2$). 
By the following proposition such a tuple is necessarily unique.

\begin{prop}\label{prop:ext primitive}
Let $L_0$ be a finite subscaled lattice,
$L$ an extension generated over $L_0$ by a primitive 
tuple $(x_1,x_2)$, and let $g=g(x_1,L_0)$.

Then $L$ is exactly the upper 
semi-lattice generated over $L_0$ 
by $x_1,x_2$. It is a finite subscaled lattice and 
one of the following happens: 
\begin{enumerate}
\item
\label{item:ext prim 1 gene}
$x_1=x_2$, $\scdim x_1<\scdim g$ and 
$\cI(L)=\cI(L_0)\cup\{x_1\}$.
\item
\label{item:ext prim 2 gene}%
$x_1\neq x_2$, $\scdim x_1=\scdim g$ and 
$\cI(L)=\left(\cI(L_0)\setminus\{g\}\right)\cup\{x_1,x_2\}$.
\end{enumerate}
\end{prop}

\begin{proof}
Let $L_1$ be the upper semi-lattice generated over $L_0$ 
by $x_1,x_2$. In order to show that $L_1=L$ it is sufficient 
to check that $L_1$ is an $\lsc$\--substructure of $L$. 

Let $p=\scdim g$ and  $q=\scdim x_1=\scdim x_2$. 
By \tcref{TC:a=(a inter b) union (a-b)}, $g^-\ll g$
hence $g^-\meet x_i\ll g$ since $g^-\meet x_i\in L_0$ 
by assumption. So $\scdim g^-<p$. We will need the following 
facts, for every $a\in L_0$: 
\begin{equation}
\label{eq:a inter xi}
g\not\leq a\quad \Ldonc\quad 
a\meet x_i\in L_0 \hbox{ \ and \ }a\meet x_i\ll x_i
\end{equation}
Indeed $g\meet a=g^-\meet a$ hence 
$a\meet x_i=(a\meet g)\meet x_i=a\meet (g^-\meet x_i)$. 
By assumption $x_i\meet g^-\in L_0$ hence $a\meet x_i\in L_0$ 
is proved. Moreover $x_i\meet (g^-\meet a)\leq x_i\meet g^-$ 
hence it suffices to check that $g^-\meet x_i\ll x_i$. 
If $x_1=x_2$ this is an assumption. Otherwise 
$x_1=g-x_2$ and $x_2=g-x_1$ are $p$\--sc\--pure
by \scref{SC:Ck(a) est k pure}. Then $x_i-g^-=x_i$ by 
\scref{SC:Ck(a) - b = Ck(a)}, since $\scdim g^-<p$, so 
$g^-\meet x_i\ll x_i$.

$L_1$ is a sublattice of $L$ then follows easily 
from (\ref{eq:a inter xi}), the distributivity law, and 
the fact that by assumption $x_1\meet x_2\in L_0\cup\{x_1\}$. 
\smallskip

Since $L_0$ is finite and $L_1$ is generated by 
$L_0\cup\{x_1,x_2\}$ as an upper semi-lattice, it follows 
immediatly that $L_1$ is finite and:
\begin{equation}
\label{eq:irred L et L0: inclusion}
\cI(L_1)\subseteq \cI(L_0)\cup\{x_1,x_2\}
\end{equation}
So for any $b'<b\in cI(L_1)$ it is easily seen that 
$\scdim(b')<\scdim b$ in $L$. 
Corollary~\ref{cor:CNS ss treillis ss LSC structure}
then shows that $L_1$ is an $\lsc$\--substructure of $L$, 
hence $L_1=L$.
\smallskip

We turn now to the description of $\cI(L)$.
If $x_1\neq x_2$ then of course $g=x_1\join x_2\nin\cI(L)$. 
Conversely if $x_1=x_2$ then (\ref{eq:irred L et L0: inclusion})
implies that $g$ is $\join$\--irreducible in $L$, so:
\begin{equation}
\label{eq:irred L et L0: g}
g\in\cI(L)\Lssi x_1\neq x_2
\end{equation}
Assume that $\cI(L_0)\not\subseteq \cI(L)$. 
Let $b\in\cI(L_0)\setminus\cI(L)$ 
and let $y_1,\dots,y_r$ ($r\geq 2$) be its $\join$\--irreducible 
components in $L_1$. By (\ref{eq:irred L et L0: inclusion}), 
each $y_i$ either belongs to $L_0$ or to $\{x_1,x_2\}$, 
and at least one of them does not belong to $L_0$. 
We may assume without loss of generality that $y_1=x_1$. 
Then $x_1\leq b$ hence $g\leq b$. If $g<b$ then $g\ll b$ 
since $b\in\cI(L_0)$ so $b-g=b$, but then we have a contradiction: 
\[y_1\leq b-g\leq b-x_1=\jjoin_{i=2}^r y_i\] 
So $b=g$. We have proved that: 
\begin{equation}\label{eq:irred L et L0: deuxieme inclusion}
\cI(L_0)\setminus\{g\}\subseteq \cI(L)
\end{equation}
The conclusion follows by combining 
(\ref{eq:irred L et L0: inclusion}), 
(\ref{eq:irred L et L0: g}),
(\ref{eq:irred L et L0: deuxieme inclusion}) 
and an obvious argument of cardinality. 
\end{proof}

It is not difficult to deduce from 
Proposition~\ref{prop:ext primitive} that any primitively 
generated extension $L$ of a finite subscaled lattice $L_0$ is 
minimal, in the sense that there is no intermediate proper 
extension $L_0\subset L_1\subset L$. Only the converse, which we 
prove now, is actually needed in the remaining of this paper.

\begin{prop}
\label{prop:extension tf de L0 fini}
Any finitely generated proper extension $L$ of a finite
$d$\--subscaled lattice $L_0$ is the union of a finite 
chain of primitively generated extensions of $L_0$. 
\end{prop}

\begin{proof}
Since $L$ is finite by the local finiteness 
Theorem~\ref{thm:TCS-treillis type fini}, it suffices to show 
that $L$ contains a primitive extension of $L_0$. Take any 
element $x$ minimal in $\cI(L)\setminus L_0$. Observe that 
if $y$ is any element of $L$ stricly smaller than $x$ then all 
the $\join$\--irreducible components of $y$ in $L$ actually 
belong to $L_0$, so $y\in L_0$. 

Let $g=g(x,L_0)$. For every $a\in L_0$, if $a<g$ 
then $a\not\leq x$ hence $a\meet x<x$, so $a\meet x\in L_0$. 
It follows that $g^-\meet x\in L_0$, hence $g^-<g$. 
In particular $g\in\cI(L_0)$. 

Since $x\in\cI(L)$ it is sc\--pure, and moreover 
$g^-\meet x<x$ implies that $g^-\meet x\ll x$ by 
\tcref{TC:a=(a inter b) union (a-b)}. 
So if moreover $x\ll g$ then we have proved that 
$(x,x)$ is primitive over $L_0$. 

On the other hand if $x\not\ll g$ let $x_1=x$ and $x_2=g-x$. 
Since $g\in\cI(L_0)$ it is sc\--pure, and so are $x_1$ and $x_2$, 
with the same sc\--dimension $p=\scdim g$ (indeed $g-x\neq\ZERO$ 
since $x<g$). Moreover $x_1\meet x_2\ll x_2$ by \tcref{TC:(x-y)-y}
so $\scdim x_1\meet x_2 < p$ hence $g-x_1=x_2$ and $g-x_2=x_1$ 
by \tcref{TC:a=(a inter b) union (a-b)}, \tcref{TC:(x union y)- z}
and the $p$\--sc\--purity of $x_1$ and $x_2$. Eventually 
$x_1\meet x_2<x$ imply that $x_1\meet x_2\in L_0$. 
We conclude that $(x_1,x_2)$ is primitive over $L_0$. 
\end{proof}

\section{Signatures}
\label{sec:signatures}

A triple $(g,q,H)$ will be called a {\df signature} in 
a finite subscaled lattice $L_0$ if and only if $g\in\cI(L_0)$, 
$q\leq\scdim g$ is a positive integer, 
and $H$ is a subset of two (not necessarily distinct)
elements $h_1,h_2\in L_0$ such that:
\begin{itemize}
\item 
$h_1\join h_2<g$ (hence $h_1\join h_2\ll g$ 
by \tcref{TC:a=(a inter b) union (a-b)})
\item
$\scdim(h_1\join h_2)<q$
\item 
If $q<\scdim g$ then $h_1=h_2$. 
\end{itemize}
By Proposition~\ref{prop:ext primitive} and the definition, 
any extension $L$ primitively generated over $L_0$ by $x_1,x_2$ 
determines a unique signature: 
\[
\sigma(L)=(g(x_1,L_0),\scdim x_1,\{x_1\meet g^-,x_2\meet g^-\})
\] 
which we call the {\df signature of $L$ in $L_0$}. 
It determines the extension $L|L_0$ as follows.

\begin{prop}\label{prop:signature et isom}
Two primitively generated extensions of a finite subscaled 
lattice $L_0$ are isomorphic over $L_0$ if and only if they 
have the same signature in $L_0$.
\end{prop}

\begin{proof}
Let $L$ (resp $L'$) be an extension of $L_0$ generated 
over $L_0$ by a primitive tuple $(x_1,x_2)$ (resp. $x'_1,x'_2$). 
If they are isomorphic over $L_0$ then obviously
$\sigma(L)=\sigma(L')$. Conversely assume that:
\[
\sigma(L)=\sigma(L')=(g,q,\{h_1,h_2\})
\]
Reordering if necessary we may assume that each 
$h_i = x_i\meet g^- = x'_i\meet g^-$. 
Let $\varphi(x_i)=x'_i$ for each $i\in\{1,2\}$, and $\varphi(a)=a$ 
for every $a\in\cI(L_0)$. 
By Proposition~\ref{prop:ext primitive} 
this is a well defined bijection from $\cI(L)$ to $\cI(L')$ 
which preserves the order, hence it extends uniquely to an 
isomorphism of upper semi-lattice $\varphi\colon L'\to L$. 
The $\llat$\--structure of a lattice being an extension 
by definition of its upper semi-lattice structure, 
this is an $\llat$\--isomorphism. 
Moreover by construction $\varphi$ preserves the sc\--dimensions 
of the $\join$\--irreducible elements of $L$ and $L'$, hence by 
Proposition~\ref{prop:CNS de lsc plongement} it is an 
$\lsc$\--isomorphism, whose restriction to $L_0$ is the identity. 
\end{proof}

\begin{prop}\label{prop:realisation d une signature}
Let $\sigma=(g,q,\{h_1,h_2\})$ be a signature in a finite 
subscaled lattice $L_0$. There exists a primitive extension 
$L$ of $L_0$ whose signature in $L_0$ is precisely $\sigma$. 
\end{prop}

\begin{proof}
We only treat the case when $q=\scdim g$. The case 
when $q<\scdim g$ (hence $h_1=h_2$) is similar, and left 
to the reader. 

Let $x_1,x_2$ be any two distinct elements in a 
set disjoint from $\cI(L_0)$ and let: 
\[
\cI=(\cI(L_0)\setminus\{g\})\cup\{x_1,x_2\}
\]
The order on $\cI(L_0)\setminus\{g\}$ inherited from $L_0$ 
can be extended to $\cI$ by stating that $x_1\not\leq x_2$, 
$x_2\not\leq x_1$, and for every 
$b\in\cI(L_0)\setminus\{g\}$ and every $j\in\{1,2\}$:
\begin{itemize}
\item 
$b<x_j\Lssi b\leq h_j$.
\item
$x_j<b\Lssi g\leq b$
\end{itemize}
For every $z\in \cI$ let $z\up=\{y\in\cI\tq z\leq y\}$. 
Let $L$ be the upper-semilattice generated in the power 
set $\cP(\cI)$ of $\cI$ by all the $z\up$'s with 
$z$ ranging over $\cI$. This is a sublattice of $\cP(\cI)$ 
such that $\cI(L)=\{z\up\tq z\in\cI\}$, that is: 
\[
\cI(L)=
\bigl(\{b\up\tq b\in\cI(L_0)\}\setminus\{g\up\}\bigr)
\cup\{x_1\up,x_2\up\}
\]
Let $\varphi(g)=x_1\up\cup x_2\up$, and 
$\varphi(b)=b\up$ for every $b\in\cI(L_0)\setminus\{g\}$. 
This application uniquely extends to an embedding of upper 
semi-lattice that we still denote $\varphi$ from $L_0$ to $L$, 
which is easily seen to be an $\llat$\--embedding. 
For every $z\in\cI$ let:
\[
D(z\up)=
\left\{
\begin{array}{ll}
\scdim(z) & \hbox{if }z\in\cI(L_0)\setminus\{g\},\\
\scdim(g) & \hbox{otherwise.}
\end{array}
\right.
\]
This is a strictly increasing map from $\cI(L)$ to $\NN$ hence 
it determines an $\lscd d$\--structure on $L$
as in Example~\ref{exple:Lnoeth}.  
Proposition~\ref{prop:CNS de lsc plongement} asserts that $\varphi$ 
is an $\lscd d$\--embedding. By construction $L$ is 
primitively generated over $L_0$ by $x_1\up, x_2\up$, 
with signature $\sigma$ in $L_0$. 

The case when $q<\scdim g$ (hence $h_1=h_2$) is similar, 
and left to the reader. 
\end{proof}

\section{Model-completion}
\label{sec:model-completion}

We call {\df super $d$\--scaled lattice} 
(resp. {\df super scaled lattice}) any 
$d$\--subscaled lattice (resp. subscaled lattice) $L$
which satisfy the following additionnal properties.
\begin{description}
\item[\bf Scaling:]
$L$ is a scaled lattice.
\item[\bf  Catenarity:]
\label{def catenarity}%
For every positive integers $r\leq q\leq p$ and every 
elements $c\leq a$, if $a$ has pure dimension $p$ and $c$ has pure 
dimension $r$ then there exists a $q$-pure element $b$ such 
that $c\leq b\leq a$.
\item[\bf Splitting:]
For every elements $b_1,b_2,a$, if $b_1\join b_2\ll a$ then 
there exists non-zero elements $a_1\geq b_1$ and $a_2\geq b_2$ 
such that:
\[
\left\{
\begin{array}{l}
a_1=a-a_2\\
a_2=a-a_1\\
a_1 \meet a_2 = b_1\meet b_2 
\end{array}
\right.
\]
\end{description}
All these properties are clearly axiomatizable in $\lscd d$,
using only finitely many $\A\E$\--formulas.
The name of the second one comes form the fact that 
in a subscaled lattice whose spectrum  is a noetherian 
topological space, this property is equivalent to the usual 
notion of catenarity, namely that any two maximal chains in 
$\Sp(L)$ having the same first and last element have 
the same length.

\begin{prop}\label{prop:plong dans superscaled}
Every subscaled lattice $\lsc$\--embeds in a superscaled lattice.
\end{prop}

\begin{proof}
Let $L_0$ be a finitely generated subscaled lattice 
and $a,b_1,b_2\in L_0$ such that $b_1\join b_2\ll a$. 
By Theorem~\ref{thm:TCS-treillis type fini}, $L_0$ is finite.
Let $g_1,\dots,g_r$ be the $\join$\--irreducible components 
of $a$ in $L_0$, and for every $i\leq r$ let $g_i^-$ be the unique 
predecessor of $g_i$ in $L_0$, $h_{i,1}=b_1\meet g_i^-$ 
and $h_{i,2}=b_2\meet g_i^-$. 

$\sigma_1=(g_1,\scdim g_1,\{h_{1,1},h_{1,2}\})$ is a signature 
in $L_0$. Proposition~\ref{prop:realisation d une signature} 
gives an extension  $L_1$ primitively generated over $L_0$ 
by $x_{1,1}, x_{1,2}$, with signature $\sigma_1$ in $L_0$. 
By construction $x_{1,1}\meet x_{1,2}$ belongs to $L_0$ and 
is stricly smaller than $g_1$ hence $x_{1,1}\meet x_{1,2}\leq g_1^-$
For every $i\geq 2$, 
$g_i\in\cI(L_0)\setminus\{g_1\}\subseteq\cI(L_1)$, and 
moreover $g_i^-$ is still the unique predecessor of $g_i$ 
in $L_1$. Indeed let $g_i^\dag$ be the union of every $c\in\cI(L_1)$ 
strictly smaller than $g_i$. Neither $x_{1,1}$ nor $x_{1,2}$ 
are smaller than $g_i$ because $g_1\not\leq g_i$, so
every such $c$ must belong to $\cI(L_0)$ by 
Proposition~\ref{prop:ext primitive}. It follows that 
$g_i^\dag=g_i^-$.

So we can repeat in $L_1$ the same construction applied 
to $g_2, h_{2,1},h_{2,2}$, and after $r$ steps we obtain 
a chain of extensions $(L_i)_{i\leq r}$ and non-zero elements
$x_{i,j}\in L=L_r$ such that for every $i$ ($1\leq i\leq r$):
\begin{eqnarray}
&
g_i-x_{i,1}=x_{i,2} 
\hbox{ \ and \ }
g_i-x_{i,2}=x_{i,1} 
& \label{eq:a egale g de x}\\
&
x_{i,1}\meet g_i^-=h_{i,1}
\hbox{ \ and \ } 
x_{i,2}\meet g_i^-=h_{i,2} 
& \label{eq:x inter a moins}\\
&
x_{i,1}\meet x_{i,2}\leq g_i^- 
& \label{eq:xi meet xj}
\end{eqnarray}
Let $x_1=\jjoin_{i\leq r}x_{i,1}$ and $x_2=\jjoin_{i\leq r}x_{i,2}$. 
For every $i\leq r$ and every $k\neq i$, 
$g_i\meet x_{k,1}\leq g_i\meet g_k\ll g_i$, so by 
\tcref{TC:x - (y union z)} and (\ref{eq:a egale g de x}):
\[
g_i- x_1 = g_i- \jjoin_{1\leq k\leq r}x_k = g_i-x_{i,1} = x_{i,2}
\]
So by \tcref{TC:(x union y)- z}, 
$a-x_1=\jjoin_{i\leq r}(g_i-x_1)=x_2$. 
Symmetrically $a-x_2=x_1$. 

(\ref{eq:x inter a moins})
and (\ref{eq:xi meet xj}) imply that
$x_{i,1}\meet x_{i,2}
=x_{i,1}\meet g_i^-\meet x_{i,2}
=h_{i,1}\meet h_{i,2}
$. Similarly for every $k\neq i$:
\[
x_{i,1}\meet x_{k,2}
=x_{i,1}\meet g_i\meet g_k\meet x_{k,2}
=x_{i,1}\meet g_i^-\meet g_k^-\meet x_{k,2}
=h_{i,1}\meet h_{k,2}
\]
So eventually:
\[
x_1\meet x_2
=\jjoin_{1\leq i \leq r}\mmeet_{1\leq k\leq r}x_{i,1}\meet x_{k,2}
=\jjoin_{1\leq i \leq r}\mmeet_{1\leq k\leq r}h_{i,1}\meet h_{k,2}
=b_1\meet b_2
\]

We have proved that for every finitely generated 
scaled lattice $L_0$ and every $a,b_1,b_2\in L_0$ 
such that $b_1\join b_2\ll a$, there exists an $\lsc$\--embeding 
$\varphi$ from $L_0$ to a subscaled lattice $L$ 
in which there exists non-zero elements $x_1,x_2$ such that, 
after identifying $L_0$ to its image: 
\[
\left\{
\begin{array}{l}
x_1=a-a_2\\
x_2=a-a_1\\
x_1 \meet x_2 = \varphi(b_1)\meet \varphi(b_2)
\end{array}
\right.
\]

On the other hand, Theorem~\ref{thm:TCS-treillis type fini} and 
Proposition~\ref{prop:linear repres}, show that every finitely 
generated subscaled lattice also $\lsc$\--embeds into some 
$\llin(X)$, which is a catenary scaled lattice. 
The model-theoretic compactness argument then implies 
that every subscaled lattice $\lsc$\--embeds in a superscaled lattice.
\end{proof}

\begin{prop}\label{prop:signature et extensions}
Let $L_0$ be a finite $\lsc$\--substructure of 
a super scaled lattice $\hat L$. 
Then for every signature $\sigma$ in $L_0$ there exists a 
primitive tuple $(x_1,x_2)\in\hat L$ such that $\sigma$ 
is the signature of $L_0\gen{x_1,x_2}$ in $L_0$. 
\end{prop}

\begin{proof}
Let $\sigma=(g,q,\{h_1,h_2\})$ be a signature in $L_0$, 
let $p=\scdim g$ and $r=\scdim(h_1\join h_2)$. 
Observe that since $g$ is 
$\join$\--irreducible in $L_0$ it is sc\--pure, 
and admits a unique predecessor $g^-$ in $L_0$. 
Let $y_1,y_2\in \hat L$ given by 
the splitting property applied to $g,h_1,h_2$. 
By construction $y_1\join y_2= g$, and since $g$
has pure sc\--dimension $p$ so does each $y_i$. 
Moreover: 
\[y_1\meet h_2\leq y_1\meet y_2=h_1\meet h_2\]
hence $y_1\meet (h_1\join h_2)=h_1\join(y_1\meet h_2)=h_1$. 
Since $h_1\join h_2<g$ it follows that $y_1\meet g^-=h_1\in L_0$, 
and symetrically 
$y_2\meet g^-=h_2\in L_0$. So $(y_1,y_2)$ is primitive over 
$L_0$, and the signature of $L_1=L_0\gen{y_1,y_2}$ in $L_0$ 
is $(g,p,\{h_1,h_2\})$. 

If $p=q$ then we are done. Assume now that $q<p$, hence 
$h_1=h_2$. For every $i\leq r$, the catenarity 
gives an element $x_i\in\hat L$ of pure sc\--dimension $p$ 
such that $\pc^i(h_1)\leq x_i\leq y_1$. Even in case $r=-1$ 
the catenarity gives $x_{-1}\in\hat L$ of pure sc\--dimension 
$p$ such that $x_{-1}\leq y_1$. Let $x=\jjoin_{-1\leq i\leq r}x_i$, 
by construction $x$ has pure sc\--dimension $p$ and 
$h_1\leq x\leq y_1$. Then $x\leq g$, and in view of their 
sc\--dimension $x\ll g$. Moreover: 
\[
h_1\leq x\meet g^- \leq y_1\meet g^- = h_1
\]
hence $h_1\meet x=h_1\in L_0$ and in view of the 
sc\--dimensions, $h_1\ll x$. So $(x,x)$ is a primitive 
tuple over $L_0$, of signature $\sigma$ in $L_0$.
\end{proof}

\begin{thm}\label{thm:model-completion}
The theory of super $d$\--scaled lattices is the 
model-completion of the theory of $d$-subscaled lattices 
\end{thm}

\begin{rem}
Obviously the theorem remains true by replacing 
everywhere $d$\--(sub)scaled lattices by 
(sub)scaled lattices of sc\--dimension at most $d$ 
(resp. exactly $d$).
\end{rem}

\begin{proof}
Since the asioms of super $d$\--scaled lattices are 
$\A\E$, it follows from Proposition~\ref{prop:plong dans superscaled} 
that every existentially closed $d$\--subscaled lattice 
is super $d$\--scaled. 
By classical model-theoretic arguments, it is then sufficient 
to show that: 
given a $d$\--superscaled lattice $\hat L$, a finitely 
generated $d$\--subscaled lattice $L$, and a common 
$\lscd d$\--substructure 
$L_0$ of $L$ and $\hat L$, there exists an $\lscd d$\--embedding 
of $L$ into $\hat L$ whose restriction to $L_0$ 
is the identity. 

By the local finiteness 
Theorem~\ref{thm:TCS-treillis type fini}, $L_0$ is finite. 
By Proposition~\ref{prop:extension tf de L0 fini} an 
immediate induction 
allows us to assume w.l.o.g. that $L$ is primitively generated 
over $L_0$. 
Let $\sigma$ be the signature of $L$ over $L_0$. By 
Proposition~\ref{prop:signature et extensions} there exists 
a primitively generated extension $L_1$ of $L_0$ in $\hat L$
whose signature is $\sigma$. Eventually $L$ is isomorphic 
to $L_1$ over $L_0$ by Proposition~\ref{prop:signature et isom}.
\end{proof}

The completions of the theory of super $d$\--scaled lattices 
are easy to classify. Let us say that a $d$\--subscaled lattice 
is {\df prime} if it does not contain any proper $d$\--subscaled 
lattice, or equivalently if it is generated by the empty set. Every 
prime $d$\--subscaled lattice is finite. By 
Corollary~\ref{cor:structures a n generateurs} their exists 
finitely many prime $d$\--subscaled lattices up to isomorphism. 

\begin{cor}\label{cor:completions}
The theory of super $d$\--scaled lattices containing (a copy of) 
a given prime $d$\--subscaled lattice is $\aleph_0$\--categorical, 
hence complete. It is also finitely axiomatisable, hence decidable. 
Since every completion of the theory of super $d$\--scaled lattices
is of that kind, the theory of super $d$\--scaled lattices
is decidable. 
\end{cor}

\begin{proof}
Let $L$, $L'$ be any two countable super $d$\--scaled lattices 
containing isomorphic prime $d$\--subscaled lattice $L_0$ and $L'_0$. 
Any partial ismorphism between $L$ and $L'$ (extending the 
given isomorphism between $L_0$ and $L'_0$) can be extended 
by va-et-vient, using exactly the same argument as in the 
proof of Theorem~\ref{thm:model-completion}. So the first 
statement is proved. The other ones are immediate consequences. 
\end{proof}

\section{Atomic scaled lattices}
\label{sec:atom-scal-latt}

A natural example of a super scaled lattice is not easy to find. 
Indeed if $X$ is any topological space in which points are closed, 
the points of $X$ are the atoms of $L(X)$, and the 
splitting axiom imply that a super scaled lattice has no atom. 
In this section we explore a possible solution to this problem, 
which lead us to a new conjecture in $p$-adic semi-algebraic 
geometry.
\smallskip

Let $\lasc=\lsc\cup\{\AT_k\}_{k\in\NN^*}$
with each $\AT_k$ a new unary predicate symbol. 
We call {\df sub-ASC-lattices} the $\lasc$\--structures 
whose $\lsc$\--reduct is a subscaled lattice and which 
satisfy the following lists of universal axioms:
\newcommand{\ascref}[1]{$\rm ASC_{\ref{#1}}$}
\begin{list}{$\rm\bf ASC_{\theenumi}:$}{\usecounter{enumi}}
\item\label{ASC:Atk donc pas Atl}
$\displaystyle
(\A k\neq l),\quad
\A a,\ \AT_k(a)\to\lnot\AT_l(a)
$
\item\label{ASC:Atk donc adb}
$\displaystyle
(\A k),\quad
\A a,a_0,\dots,a_{2^k},\quad \AT_k(a)\Lto$ 

\centerline{
$\displaystyle
\scdim a= 0 
\hbox{$\bigwedge$} \left[
  \Bigl(\jjoin_{0\leq i\leq 2^k}a_i=a\Bigr)\to
  \Bigl(\jjoin_{0\leq i<j\leq 2^k}a_i=a_j\Bigr) \right]
$}
\item\label{ASC:Atk et les Atl}
$\displaystyle
(\A k,n),\quad
\A a,a_1,\dots,a_n,$

\centerline{
$\displaystyle
\left[
\Bigl(a=\jjoin_{1\leq i\leq n}a_i\Bigr)
\hbox{$\bigwedge$} 
\Bigl(\mmeet_{1\leq j<i\leq n}a_i\meet a_j =\ZERO\Bigr)
\right]
\Lto\qquad\qquad
$}

\centerline{
$\qquad\qquad\biggl[
\AT_k(a)\Lliff
\hbox{$\ddisj\limits_{l_1+\cdots+l_n=k}\;\cconj\limits_{1\leq i\leq n}$}
\AT_{l_i}(a_i)
\biggr]
$}
\end{list}

For any $\lasc$\--structure $L$ we denote by 
$\AT_k(L)$ the set of elements $a$ in $L$ such that 
$L\models\AT_k(a)$, and we let 
$\AT_0(L)=L\setminus\bigcup_{k>0}\AT_k(L)$. 
If $L$ is a sub-\ASC-lattice then 
\ascref{ASC:Atk donc pas Atl} asserts that 
$(\AT_k(L))_{k\in\NN}$ is a partition of $L$. 
For any $a\in L$ we then define $\asc(a)$ as the unique 
$k\in\NN$ such that $a\in\AT_k(L)$. The other axioms 
of sub-\ASC-lattices have the following meaning:
\begin{itemize}
\item
If $\asc(a)=k>0$ then \ascref{ASC:Atk donc adb} asserts
that $L(a)$ is a boolean algebra generated by $n$ atoms 
$a_1,\dots,a_n$ with $n\leq k$, and \ascref{ASC:Atk et les Atl} 
then implies that each $\asc(a_i)>0$ and for every $b\in L(a)$:
\[
\asc(b)=\sum_{c\in\cC(b)}\asc(c)
\]
where $\cC(b)$ is the set of atoms in $L(a)$ (as well as in $L$) 
smaller than $b$.
\item
Conversely \ascref{ASC:Atk et les Atl} implies that 
if $L(a)$ is finite and $\asc(c)>0$ for each atom $c$ in $L(a)$:
\[
\asc(a)=\sum_{c\in\cI(L(a))}\asc(c)
\]
In particular if $a\neq\ZERO$ then $\asc(a)>0$.
\end{itemize}

\begin{rem}\label{rem:lasc plongement}
It follows immediatly that an $\lsc$\--embedding of 
sub-\ASC-lattices $\varphi\colon L\to L'$ is an $\lasc$\--embedding 
if and only if for every $k>0$ and every {\it atom} $a\in L$:
\[
L\models\AT_k(a)
\Lssi
L'\models\AT_k(\varphi(a))
\]
\end{rem}

Every natural example of sub-\ASC-lattice $L$ satisfy the following 
additionnal property, which imply \ascref{ASC:Atk donc pas Atl} to 
\ascref{ASC:Atk et les Atl}:
\begin{list}{$\rm\bf ASC_{\theenumi}:$}{\usecounter{enumi}}
\setcounter{enumi}{-1}
\item
\label{ASC:def ASC}
$L$ is a scaled lattices, and for every $k>0$, $\AT_k(L)$ 
is  the set of elements of $L$ which are the join of exactly 
$k$ atoms in $L$. 
\end{list}
We call {\df ASC\--lattices}\footnote{Of course 
we will show that the sub-\ASC-lattices are precisely 
the $\lasc$\--substructures of \ASC-lattices.}  
the sub-\ASC-lattices which satisfy \ascref{ASC:def ASC}. 
Every subscaled lattice $L$ admits a unique structure 
of \ASC-lattice which is an extension by definition of its lattice 
structure. We denote by $L^{\rm At}$ this expansion of $L$. 

\begin{prop}[Linear representation]\label{prop:repres lin ASC}
Let $\KK$ be an infinite field and let $L_0$ be a finite 
sub-\ASC-lattice. 
For any positive integer $N$ there exists a positive integer $m$, 
not depending on $N$, and an $\lsc^*$\--embedding 
$\varphi_N\colon L_0\to\lalin(\KK^m)$
such that for every element $a$ of $L_0$ of sc\--dimension zero:
\begin{eqnarray*}
\asc(a)>0 &\Ldonc& \asc(\varphi_N(a))=\asc(a)\\
\asc(a)=0 &\Ldonc& \asc(\varphi_N(a))\geq N
\end{eqnarray*}
\end{prop}

\begin{proof}
We prove by induction on tuples $(r,s)$ ordered 
lexicographically that the proposition holds for every finite 
sub-\ASC-lattice $L_0$ having $r$ non-zero $\join$\--irreducible 
elements, $s$ of which have the same sc\--dimension as $L_0$. 

If $r=0$ then $s=0$ and the 
unique embedding of $L_0=\{\ZERO\}$ into $\lalin(\KK^0)$ 
has the required property. So let us assume that $r\geq 1$ and that 
the result is proved for every $(r',s')<(r,s)$. 
Let $a_1\dots,a_r$ be the elements of $\cI(L_0)$ ordered 
by increasing sc\--dimension, so $d=\scdim L=\scdim a_r\geq 0$.
\smallskip

{\bf Case $d=0$.}
Then $L$ is a boolean algebra, and $a_1,\dots,a_r$ are its atoms. 
For every $i\leq r$ we choose a finite subset $A_i$ of 
$\KK\setminus\smash{\bigcup\limits_{j<i}A_j}$ such that:
\begin{itemize}
\item
If $\asc(a_i)>0$, $A_i$ has $\asc(a_i)$ elements, so 
$\asc(A_i)=\asc(a_i)$.
\item
If $\asc(a_i)=0$, $A_i$ has $N$ elements, so $\asc(A_i)=N$. 
\end{itemize}
Clearly $m=1$ does not depend on $N$, and the map 
$\varphi$ which maps each $a_i$ to $A_i$ extends uniquely 
to an $\lsc^*$\--embedding of $L_0$ into $\lalin(\KK)$ which 
has the required properties.
\smallskip

{\bf Case $d>0$.}
The upper semi-lattice $L_0^-$ generated by $a_1,\dots,a_{r-1}$
is an $\lasc^*$\--substructure of $L_0$ to which the induction 
hypothesis applies. This gives a positive integer $m$ not depending 
on $N$ and an $\lsc^*$\--embedding $\varphi\colon L_0^-\to\lalin(\KK^m)$
having the required properties. One can extend $\varphi$ to an 
embedding $\bar\varphi$ of $L_0$ into some $\lalin(\KK^{m+p})$ 
exactly like in the proof of Proposition~\ref{prop:linear repres}. 
The integer $m+p$ does not depend on $N$ and $\bar\varphi$ 
inherits from $\varphi$ the required properties because 
all the elements of $L$ with sc\--dimension zero belong 
to $L_0^-$. 
\end{proof}

Let $\sazar(\KK,d)$ (resp. $\salin(\KK,d)$, 
resp. $\sadef(\KK,d)$) the class of all \ASC-lattices $L^{\rm At}$ 
with $L$ ranging over $\szar(\KK,d)$, (resp. $\slin(\KK,d)$, 
resp. $\sdef(\KK,d)$). 

\begin{cor}\label{cor:repres ASC non standard}
For any infinite field $\KK$ and positive integer $d$, 
the universal theories of $\sazar(\KK,d)$ and $\salin(\KK,d)$)
are exactly the theory of sub-\ASC-lattices. 

The same holds for $\sadef(\KK,d)$ if moreover $\KK$ is a henselian 
valued field of characteristic zero, a real closed field or an 
algebraically closed field. 
\end{cor}

\begin{proof}
Since $\salin(\KK,d)$ is contained in the other classes, all of 
which are contained in the class of \ASC-lattices, it suffices 
to prove that every sub-\ASC-lattice $\lasc$\--embeds into 
an ultraproduct of elements of $\salin(\KK,d)$. 

For any positive integer $N$ let $\varphi_N:L_0\to\lalin(\KK^m)$ 
be the $\lsc^*$\--embedding given by 
Proposition~\ref{prop:repres lin ASC}. 
We still denote by $\varphi_N$ the induced $\lsc$\--embedding 
of $L_0$ into $\lalin(X_N)$ where $X_N=\varphi_N(\UN_{L_0})$. 
Let $\cU$ be a non principal ultrafilter 
in the boolean algebra of subsets of $\NN$, and 
$L=\prod_{N\in\NN}\lalin(X_N)/\cU$. 
Then $\varphi=\prod_{N\in\NN}\varphi_N/\cU$ is 
an $\lsc$\--embedding of $L_0$ into the ultraproduct $L$. 
Let $a$ be an atom of $L_0$ and $k=\asc(a)$. 

If $k>0$ then for every $N\geq k$, 
$\lalin(X_N)\models\AT_k(\varphi_N(a))$ by construction. So 
$L\models\AT_k(\varphi(a))$, that is $\asc(\varphi(a))=k$. 

If $k=0$, let $l$ be any stricly positive integer. For every 
$N\geq l$, $\lalin(X_N)\models\AT_N(\varphi_N(a))$ 
by construction, hence $lalin(X_N)\not\models\AT_l(\varphi_N(a))$. 
So $L\not\models\AT_l(\varphi(a))$, and 
this being true for every $l>0$ it follows that $\asc(\varphi(a))=0$. 

By Remark~\ref{rem:lasc plongement}, $\varphi$ is then an 
$\lasc$\--embedding.
\end{proof}

Let us call {\df super ASC\--lattices} those \ASC-lattices 
which satisfy the following axioms:
\begin{description}
\item[\bf Atomicity:] 
Every element $x$ is the least upper bound in $L$ of the set of 
atoms of $L$ smaller that $x$. 
\item[\bf Catenarity:]
For every positive integers $r\leq q\leq p$ and every 
$c\leq a$, if $a$ has pure dimension $p$ and $c$ has pure 
dimension $r$ then there exists a $q$-pure element $b$ such 
that $c\leq b\leq a$.
\item[\bf ASC-Splitting:]
For every $b_1,b_2,a$, if $b_1\join b_2\ll a$ and $\pc^0(a)=\ZERO$ 
there exists non-zero elements $a_1\geq b_1$ and $a_2\geq b_2$ 
such that:
\[
\left\{
\begin{array}{l}
a_1=a-a_2\\
a_2=a-a_1\\
a_1 \meet a_2 = b_1\meet b_2 
\end{array}
\right.
\]
\end{description}
The class of super-\ASC-lattices of sc\--dimension at most $d$ 
(resp. exactly $d$) is clearly axiomatisable by 
$\A\E$\--formulas in $\lasc$. We are going to show that its 
theory is the model-completion of the theory of sub-\ASC-lattices 
of dimension at most $d$ (resp. exactly $d$). 

\begin{rem}\label{rem:infinite d'atomes}
An immediate consequence of the atomicity axiom
is that for every elements $x,y$ in a super \ASC-lattice $L$ 
such that $y\ll x$, there are infinitely many atoms $a\in L$ such 
that $a\leq x$ and $a\meet y=\ZERO$. 
Indeed $y<x$ hence by the atomicity axiom there 
is an atom $a_1\in L$ such that $a_1\leq x$ and $a_1\not\leq y$. 
Then $a_1\leq x$ and $a_1\meet y<a_1$ hence 
$a_1\meet y=\ZERO$ (because $a_1$ is an atom).
Moreover $\scdim y\join a_1<\scdim x$ because $y\ll x$, 
hence $y\join a_1<x$ so the same argument applies 
to $x$ and $y\join a_1$. It gives another atom $a_2\leq x$ 
such that $a_2\meet y=\ZERO$, and so on. 
\end{rem}

Primitive tuples and primitive extensions are defined for 
sub-\ASC-lattices exactly like for subscaled lattices. 
If $L$ is an extension of a sub-\ASC-lattice $L_0$ and 
$x\in L$ we denote now $L_0\gen x$ the $\lasc$\--substructure 
of $L$ generated by $L_0\join\{x\}$ (that is the subscaled 
lattice generated by $L_0\cup\{x\}$ endowed with the 
$\lasc$\--structure induced by $L$).

We define {\df ASC\--signatures} in a finite sub-\ASC-lattice $L_0$ 
as triples $(g,p,H)$ with $H$ a set of non-necessarily distinct 
ordered pairs $(h_1,k_1)$, $(h_2,k_2)$ such that $(g,q,\{h_1,h_2\})$ 
is a signature in the $\lsc$\--reduct of $L_0$, $k_1,k_2$ are 
positive integers and: 
\begin{enumerate}
\item
If $q<\scdim g$ then $(h_1,k_1)=(h_2,k_2)$.
\item
If $q\neq 0$ then $k_1=k_2=0$.
\item 
If $k_1\neq 0$, $k_2\neq 0$ and $\scdim g=0$ 
then $\asc(g)=k_1+k_2$
\item
If $k_1=0$ or $k_2=0$ then $\asc(g)=0$.
\end{enumerate}

\begin{exmp}
Let $L_0$ be a finite sub-\ASC-lattice, and $L$ an extension 
of $L_0$ generated by a primitive tuple $(x_1,x_2)$. Let 
$(g,p,\{h_1,h_2\})$ be the signature of $L$ in $L_0$ 
(in the sense of subscaled lattices) and $k_i=\asc(x_i)$. 
Then $(g,p,\{(h_1,k_1),(h_2,k_2)\}$ is easily seen to be 
an \ASC-signature in $L_0$, uniquely determined by $L$.
We call it the {\df ASC\--signature of $L$ in $L_0$}. 
\end{exmp} 

The same argument as in Proposition~\ref{prop:signature et isom} 
shows (using Remark~\ref{rem:lasc plongement} in additition to 
Proposition~\ref{prop:CNS de lsc plongement}) that two primitively 
generated extensions of a finite sub-\ASC-lattice $L_0$ are 
isomorphic over $L_0$ if and only if they have the 
same signature in $L_0$. 

\begin{prop}\label{prop:signature et extensions ASC}
Let $L_0$ be a finite $\lasc$\--substructure of 
an $\aleph_0$\--saturated super \ASC-lattice $\hat L$. 
Then for every \ASC-signature $\sigma_{\rm At}$ in $L_0$ 
there exists a primitive tuple $(x_1,x_2)\in\hat L$ 
such that $\sigma_{\rm At}$ is the signature of 
$L_0\gen{x_1,x_2}$ in $L_0$. 
\end{prop}

\begin{proof}
Let $\sigma_{\rm At}=(g,q,\{(h_1,k_1),(h_2,k_2)\})$ be 
an \ASC-signature in $L_0$, let $\sigma=(g,q,\{h_1,h_2\})$ and 
$p=\scdim g$. In each of the following cases, the 
fact that the constructed tuple $(x_1,x_2)$ is primitive 
over $L_0$ and generates an extension whose signature 
in $L_0$ is precisely $\sigma_{\rm At}$ is straighforward.
\medskip

{\bf Case $p\geq 1$.} 
Then $g$ is not an atom, so the \ASC-splitting property 
applies to $g, h_1, h_2$. Moreover by 
Remark~\ref{rem:infinite d'atomes}, there are infinitely 
many atoms $a$ in $L$ smaller than $g$ and not in $L_0$. 
By $\aleph_0$\--saturation it follows that there exists 
$x_0\in L\setminus L_0$ such that $\scdim x_0=0$ and 
$L\not\models\AT_k(a)$ for every $k>0$, that is $\asc(x_0)=0$.
\smallskip

If $q\geq 1$ then the construction in the proof of 
Proposition~\ref{prop:signature et extensions} applies here. It
gives a primitive tuple $(x_1,x_2)$ in $L$ such that the 
$\lsc$\--substructure $L_1$ of $L$ generated by $(x_1,x_2)$ 
over $L_0$ has signature $\sigma$ in $L_0$. Moreover each 
$k_i=0$ since $p\neq 0$, and on the other hand each $\asc x_i=0$ 
since $\scdim x_i=p\geq 1$.
\smallskip

Otherwise $q=0<p$ hence $(h_1,k_1)=(h_2,k_2)$. If $k_1\neq 0$ 
let $x_1=x_2=a_1\join\dots\join a_{k_1}$ with the $a_i$'s 
being any distinct atoms in $L$ smaller than $g$ and not 
belonging to $L_0$. If $k_1=0$ let $x_1=x_2=x_0$. 
\medskip

{\bf Case $p=0$.}
Then $q=0$ and $h_1=h_2=\ZERO$. 
\smallskip

If $k_1$ and $k_2$ are non-zero then $\asc(g)=k_1+k_2$ 
hence $L(g)$ contains $k_1+k_2$ atoms. Let $x_1$ be the join 
of $k_1$ of them, and $x_2$ be the join of the others. 
\smallskip

Eventually if $k_1=0$ or $k_2=0$ then $\asc(g)=0$ so $L(g)$ contains 
infinitely many atoms, hence by saturation there exists 
an element $x$ in $L$ smaller than $g$ such that both 
$x$ and $g-x$ are non-zero and $\asc(x)=\asc(g-x)=0$. 
If $k_1=k_2=0$  let $(x_1,x_2)=(x,g-x)$. 
If $k_1\neq 0$ let $x_1=a_1\join\dots\join a_{k_1}$ with 
the $a_i$'s atoms in $L$ smaller than $g$ and $x_2=g-x_1$. 
If $k_2\neq0$ exchange $k_1$ and $k_2$. 
\end{proof}

\begin{prop}
Every sub-\ASC-lattice embeds in a super-\ASC-lattice.
\end{prop}

\begin{proof}
Obviously every finitely generated substructure of a
sub-\ASC-lattices is finite by the local finiteness 
Theorem~\ref{thm:TCS-treillis type fini} because $\lsc$ 
and $\lasc$ have the same function symbols. 

Let $L_0$ be a finitely generated sub-\ASC-lattice, 
and $a,b_1,b_2\in L_0$ such that $b_1\join b_2\ll a$ 
and $\pc^0(a)=\ZERO$. So $L_0$ is finite, let $g_1,\dots,g_r$ 
be the $\join$\--irreducible components of $a$ in $L_0$. 
As in the proof of Proposition~\ref{prop:plong dans superscaled} 
we construct a tower of $\lsc$\--extensions $(L_i)_{i\leq r}$ 
so that each $L_r$ contains elements $x_1,x_2$ satisfying 
\[
\left\{
\begin{array}{l}
x_1=a-a_2\\
x_2=a-a_1\\
x_1 \meet x_2 = \varphi(b_1)\meet \varphi(b_2)
\end{array}
\right.
\]
and each $L_i$ is generated over $L_{i-1}$ by a primitive 
tuple of elements $x_{i,1}\neq x_{i,2}$ such that 
$g_i=g(x_{i,1},L_{i-1})=g(x_{i,2},L_{i-1})$. 
Proposition~\ref{prop:ext primitive} then implies that: 
\[
\cI(L_r)=\bigl(\cI(L_0)\setminus\{g_1,\dots,g_r\}\bigr)
\cup\{x_{i,1},x_{i,2}\}_{1\leq i\leq r}
\]
By assumption all the $g_i$'s have non-zero sc\--dimension, 
hence the $x_{i,1}$'s and $x_{i,2}$'s do the same. In particular 
$L_0$ and $L_r$ have the same elements of sc\--dimension zero. 
Defining $\AT_k(L_r)=\AT_k(L_0)$ for every $k$ then endows 
$L_r$ with an $\lasc$\--structure which makes it a sub-\ASC-lattice 
and an $\lasc$\--extension of $L_0$. 

On the other hand every $\lalin(X)$ is an \ASC-lattice 
satisfying the atomicity and catenarity properties, 
hence Corollary~\ref{cor:repres ASC non standard} implies 
that every finitely generated sub-\ASC-lattice embeds 
in an \ASC-lattice having these two properties. 

The conclusion follows by the model-theoretic compactness 
argument. 
\end{proof}

\begin{thm}\label{thm:model-completion ASC}
The theory of super \ASC-lattices of sc\--dimension at most $d$ 
(resp. exactly $d$) is the model-completion of the theory of 
\ASC-lattices of dimension at most $d$ (resp. exactly $d$).
It admits $\aleph_0$ completions, each of which is decidable, 
and it is decidable.
\end{thm}

\begin{proof}
The proof of the first statement is very similar to 
Theorem~\ref{thm:model-completion} and 
Corollary~\ref{cor:completions}, with the only difference that 
in the embedding argument we have to assume the 
super \ASC-lattice $\hat L$ to be $\aleph_0$\--saturated, 
in order to apply Proposition~\ref{prop:signature et extensions ASC} 
in place of Proposition~\ref{prop:signature et extensions}.

The last statement then follows from the remark that 
there are $\aleph_0$ prime sub-\ASC-lattices of dimension 
at most $d$ (resp. exactly $d$). Indeed there are 
finitely many subscaled lattices of dimension 
at most $d$ (resp. exactly $d$) and on each of them $\aleph_0$ 
different structures of sub-\ASC-lattices. 
\end{proof}

Let $L$ be any super \ASC-lattice and $L'$ the quotient of $L$ 
by the equivalence relation $x\sim y$ if and only if 
$x-y$ and $y-x$ are the join of finitely many atoms. 
Then it is easily seen that $L'$ is a super scaled lattice. 
So the problem of finding a natural example of super scaled lattice 
boyls down to finding a natural example of super \ASC-lattice. 
This and related ideas lead us to the following conjecture:

\begin{conj}\label{conj:splitting p-adique ameliore}
Let $K$ be a $p$\--adically closed field and $A$ be an infinite 
definable subset of $K^n$ which is open in its closure. 
Let $(B_k)_{k\leq q}$ be a finite collection of closed definable 
subsets of $\bar A\setminus A$. Then there exists a collection 
$(A_k)_{k\leq q}$ of non-empty definable subsets of $A$ clopen in $A$ 
such that: 
$$\A k\leq q,\quad \overline{A_k}=A_k\cup B_k$$
\end{conj}

It is an elementary exercise to deduce from this conjecture 
that $\ladef(X)$ models the ASC\--splitting property for 
every definable set $X$ over $\KK$. 
Moreover the $\lsc$\--substructure of $\ladef(\KK^n)$ generated 
by the empty set is simply the two element lattice with the 
obvious $\lsc$\--structure (because $\KK^d$ is $d$\--pure). 
This gives an explicit recursive axiomatisation of $\ladef(\KK^n)$. 

\begin{cor}[Modulo Conjecture~\ref{conj:splitting p-adique ameliore}]
Let $\KK$ be a $p$\--adically closed field, then $\ladef(\KK^n)$ 
is a super-ASC-lattice. In particular its complete theory is 
decidable and eliminates the quantifier in $\lasc$. 
\end{cor}


\begin{thebibliography}{vdD89}

\bibitem[Dar04]{prep-latt-2004}
Luck Darni{\`e}re.
\newblock Model-completion of scaled lattices.
\newblock {\em Pr\'e\-pu\-bli\-ca\-tion Ma\-th\'e\-ma\-tiques d'An\-gers 191.} May 2004.
\newblock {\tt http://math.univ-angers.fr}.

\bibitem[Grz51]{grze-1951}
Andrzej Grzegorczyk.
\newblock Undecidability of some topological theories.
\newblock {\em Fund. Math.}, 38:137--152, 1951.

\bibitem[Joh82]{John-1982}
Peter~T. Johnstone.
\newblock {\em Stone spaces}, volume~3 of {\em Cambridge Studies in Advanced
  Mathematics}.
\newblock Cambridge University Press, Cambridge, 1982.

\bibitem[vdD89]{Drie-1989}
Lou van~den Dries.
\newblock Dimension of definable sets, algebraic boundedness and {H}enselian
  fields.
\newblock {\em Ann. Pure Appl. Logic}, 45(2):189--209, 1989.
\newblock Stability in model theory, II (Trento, 1987).

\end{thebibliography}

\end{document}